\documentclass[final,1p,times]{elsarticle}

\usepackage[super,compress]{cite}
\usepackage{url}
\usepackage{hyperref}
\usepackage{verbatim}
\usepackage{graphicx}%
\usepackage{amsmath,amssymb,amsfonts}%
\usepackage{amsthm}%
\usepackage{mathrsfs}%
\usepackage[title]{appendix}%
\usepackage{xcolor}%
\usepackage{todonotes}
\usepackage{subcaption}
\usepackage{color}
\usepackage{tikz}
\usetikzlibrary{automata,positioning,arrows}
\usepackage{longtable}
\usepackage{ulem}
\usepackage{subcaption}
\usepackage{lscape}

\newtheorem{theorem}{Theorem}[section]
\newtheorem{remark}[theorem]{Remark}

\usepackage{array}
\newcounter{rowno}
\allowdisplaybreaks

\journal{Chaos, Solitons \& Fractals}

\begin{document}
	
	\begin{frontmatter}
		
		
		
		\title{Contribution of the 2021 COVID-19 Vaccination Regime to COVID-19 Transmission and Control in South Africa: \textit{A Mathematical Modeling Perspective.}}
		
		\author[UP]{Tesfalem A. Tegegn}
		\affiliation[UP]{organization={Department of mathematics and applied mathematics\\University of Pretoria},
			addressline={Private bag X20 Hatfield},
			city={Pretoria},
			postcode={0028},
			state={Gauteng},
			country={South Africa}}
		
		\author[UFS]{Yibeltal A. Terefe} 
		
		\affiliation[UFS]{organization={Department of Mathematics and Applied Mathematics \\University of the Free State},
			addressline={P.O. Box 339}, 
			city={Bloemfontein},
			postcode={9300}, 
			state={Free State},
			country={South Africa}}
		
		\begin{abstract}
			This study assesses the impact of COVID-19 vaccines through an epidemiological model that quantifies their role in pandemic control. We analyze an SIR-based model incorporating vaccination effects and calculate the basic reproduction number, $\mathcal{R}_0$. Our findings indicate that
			under conditions of imperfect vaccination and non-permanent immunity from recovery, a backward bifurcation may occur when $\mathcal{R}_0<1$. Conversely, with full immunity from vaccination and lasting immunity post-recovery, the disease-free equilibrium is globally asymptotically stable for $\mathcal{R}_0<1$. Numerical simulations support these theoretical results. The model is calibrated using South African data from Johns Hopkins University, covering the period from February 17 to August 5, 2021. Results show that vaccine effectiveness in preventing infection was below 50\%, consistent with the CDC’s February 2024 report indicating an improved infection-protection rate of 54\% for newly produced vaccines. Additionally, our findings demonstrate that vaccines significantly enhanced recovery rates and reduced both mortality and recovery time, aligning with CDC reports. A sensitivity analysis highlights key parameters affecting $\mathcal{R}_0$, offering insights for
			policymakers on optimizing vaccination strategies.  
		\end{abstract}
		
		
		
		
	\end{frontmatter}


	\section{Introduction} \label{intro}
	
	COVID-19 is an infectious disease caused by SARS-CoV-2 virus, which  prompted the World Health Organization (WHO) to declare Public Health Emergency of International Concern on 30 January 2020, and  latter a pandemic on 11 March 2020. 
	   First detected in December 2019 in Wuhan, China, the virus has since
	  spread   to all  seven continents and every country recognized by the United Nations as a sovereign state. According to the \href{https://www.worldometers.info/coronavirus/}{worldometers database}\cite{worldometer},    COVID-19 infected well over  700 million  people and claimed   over 7 million lives  as of 24 July 2024.  %
	
	Measures implemented to curb the spread of COVID-19 including large-scale lockdowns,
	social distancing and others,  significantly disrupted family and social life   while badly affecting     the global economy by crippling global and local supply chains \cite{Atalan2020lockdown,ceylan2020historical,chaudhary2020effect, Das2022socioeconomic, deb2022economic, mou2020research, singh2021health}. %
 	Global financial institutions,  including the World Bank and International Monetary Fund (IMF),   repeatedly revised their projections for  global  economic growth during the pandemic.
	A  \href{https://www.mckinsey.com/business-functions/strategy-and-corporate-finance/our-insights/the-coronavirus-effect-on-global-economic-sentiment}{McKinsey \& Company} report (October 2021)\cite{mckinsey} highlighted the  COVID-19 related  concerns, particularly travel restrictions and supply chain disruption, were perceived by both corporate leaders and government officials as the most significant threats to domestic and corporate economic growth.  

	 The SARS-CoV-2 virus primarily spreads from person to person through small respiratory droplets, known as aerosols, which are released when an infected individual coughs, sneezes, talks, or engages in similar activities \cite{who}.
		 Initially, due to the relatively large size of the virus, it was assumed that transmission could occur only via larger droplets, which were believed to remain airborne for only a short time and to travel no more than two meters from the source.
		  However, Van Doremalen et al. \cite{van2020aerosol} demonstrated that SARS-CoV-2 can remain airborne for at least three hours, with a viral load sufficient to cause infection. 
		  This finding indicates that the virus is indeed airborne, helping to explain its rapid spread and the limited effectiveness of control measures such as lockdowns and social distancing. Another significant mode of transmission is through contact with contaminated surfaces, such as doorknobs, stair railings, and elevator buttons, followed by touching the nose, mouth, or eyes before hand hygiene. Studies have shown that the virus can remain viable on hard surfaces for varying durations depending on the material \cite{van2020aerosol}.

	Since the end of the $18^{\rm th}$ century, mathematical models have been in use in public health areas and have  provided useful information for policy makers. 
	These models also assist public health professionals by generating valuable information to design mitigation and control strategies, see \cite{kramer2010modern}. 
In the context of COVID-19 pandemic, several  mathematical models have been proposed and analysed to understand the transmission dynamics of the disease. 
	One of the earliest works in this regard is \cite{wu2020nowcasting} in which the  authors used   a simple SEIR model to estimate the pandemics trajectory. 
	Subsequent studies, such as  \cite{Garba2020} proposed a modified SIR   model to analyze the role of environmental contamination by infected individuals. In \cite{kassa2020}, an SIR-based model was proposed to assess the effectiveness of mitigation strategies    recommended by the WHO. Similarly,     \cite{kassa2023}  proposed  SEIR-based model   to analyze the impact  of cross-boarder migration, with and without screening at boarder crossings, on disease prevalence  in the host community.  %
	 Indeed, in the context of pandemics, one of the key contributions of mathematical research is to project the epidemic trajectory and recommend effective control mechanisms to support public health policy decisions.

	Vaccines have long played a crucial role in controlling pandemics such as smallpox, yellow fever, and influenza. Although vaccines do not offer 100\% protection, they significantly reduce the risk of infection, severe illness, and death. In the case of COVID-19, it is well established that the efficacy\footnote{According to the \href{https://www.who.int/news-room/feature-stories/detail/vaccine-efficacy-effectiveness-and-protection}{World Health Organization} \cite{who}, vaccine efficacy refers to how well a vaccine prevents disease under ideal, controlled clinical trial conditions, whereas vaccine effectiveness describes how well it performs in real-world settings.} of vaccines is high for a limited period. They have been shown to be particularly effective in preventing severe illness, reducing the risk of infection, and lowering mortality associated with the virus. Despite not being perfect, widespread vaccination campaigns were carried out globally, and they played a critical role in bringing the COVID-19 pandemic under control.

	The primary objective of this study is not to nowcast or forecast the trajectory of the COVID-19 outbreak, but rather to assess the contribution of  vaccines administered in South Africa  during the period 17 March to 5 August 2021. 
	The specific time frame   was chosen  solely based on  data availability; a detailed explanation is provided in subsection \ref{sec: about the data}. 
	To achieve this goal, we propose an SIRS-based epidemiological model  
	in which the susceptible population is divided into vaccinated and unvaccinated classes based on their vulnerability. 
	A vaccinated individual is defined as someone   who has received   at least one dose  of any   vaccine  legally approved and administered in the country. 
	The model accounts for important real-world considerations, including the fact that  vaccines do not offer complete protection  and reinfection is possible.   
	The model also considers separate compartments for infectious people from the two susceptible groups so as to be able to analyse how vaccines helped the infected to recover.   
	Additionally, separate compartments are included for infectious individuals originating from both susceptible groups, enabling a detailed analysis of how vaccination influences recovery outcomes.
	Overall, this study introduces a novel epidemic model that classifies infected individuals into vaccinated and unvaccinated subgroups, further stratified by symptomatic and asymptomatic status. This structure enables a targeted evaluation of vaccine efficacy in facilitating recovery, an aspect that is not commonly addressed in existing literature. The full model formulation and associated assumptions are presented in Section \ref{modelform}.

	   Section \ref{analysis} is dedicated to the qualitative and quantitative analysis of the proposed model. In this section, we compute the basic reproduction number, identify the disease-free and endemic equilibrium points, and analyze their stability, among other aspects. For the numerical analysis, we used Python libraries such as Pandas, NumPy, lmfit, and others to download data from \cite{JH_data}, preprocess it, fit it to the model, and perform local sensitivity analysis. A comprehensive discussion on the data, primary parameter estimation, model fitting, sensitivity analysis, and interpretation of the results is presented in Sections \ref{fitting}, \ref{sensitivity}, and \ref{conclusion}, respectively.

	\section{Model Formulation} \label{modelform}

	In this section, we propose a mathematical model designed to capture the contribution of vaccines  in two key aspects: reducing the risk of infection, and, in cases where infection occurs despite vaccination, mitigating symptoms and enhancing the recovery rate.
	
	To formulate the model, we begin by  assuming that the population is mixed homogeneously. 
	Based on   individuals'  vulnerability to infection and infection status,   the total population  is divided into eight mutually exclusive  compartments.
	Susceptible individuals who have not received COVID-19  vaccines are placed in $S$ class. Those who have received at least one dose of any approved COVID-19 vaccine are assigned to the $V$ class. 
	Unvaccinated individuals who become infected are classified as either symptomatic  ($I$) or asymptomatic ($A$). 
	In cases where vaccination does not prevent infection, vaccinated individuals may still become infected and are classified as asymptomatic ($A_1$)  or symptomatic ($I_1$).  
	Infectious individuals who are isolated, whether in hospitals or temporary treatment centers or own facilities, are placed in the $Q$ class. 
	Finally, individuals who have recovered from COVID-19 are assigned to the $R$ class. 
	
	Thus, the total population at time $t$,  denoted by $N(t)$,   is given by 
	\begin{eqnarray}
		N(t)=S(t) +V(t) +A(t) +A_1(t) +I(t) +I_1(t) +Q(t) +R(t) , \label{eq: total population}
	\end{eqnarray}
	%
	where $S(t),$ $V(t),$ $A(t)$ $A_1(t) ,$ $I(t) ,$ $I_1(t) ,$ $Q(t) $ and $R(t)$ represent the number of individuals in compartments $S,$ $V,$ $A ,$ $A_1,$ $I,$ $I_1,$ $Q $ and $R$, respectively, at  time $t$.  For convenience, we will hereafter refer to these variables as $S,$ $V,$ $A,$ $A_1,$ $I,$ $I_1,$ $Q,$ and $R$ instead of explicitly writing them as functions of time, unless otherwise necessary.
	
	In formulating the model, we excluded both the exposed (latent) compartment and the contribution of environmental contamination in order to facilitate a more tractable mathematical analysis.   
	More specifically, although COVID-19 is known to transmit through both direct human-to-human contact and indirect contact via contaminated surfaces, this model considers only direct transmission.  
That is, we focus on transmissions arising from interactions between susceptible individuals and infectious individuals, such as through inhalation of virus laden aerosols or physical contacts, such as handshake.
	
	Given that the model includes four classes of infectious individuals, namely, symptomatic unvaccinated ($I$), asymptomatic unvaccinated ($A$), symptomatic vaccinated ($I_1$), and asymptomatic vaccinated ($A_1$)  who may interact with the susceptible individuals, 
	we define the force of infection as:  
	\begin{eqnarray}\label{FF}
		\lambda =\beta \frac{I +\nu A +\nu_1A_1 +\kappa I_1 }{N },
	\end{eqnarray}
	where $\beta$ denotes the transmission rate, and $\nu,\;\nu_1$ and $\kappa$ are modification parameters that account for the relative infectiousness of compartments $A,\,A_1,$ and $I_1,$ respectively, compared to the symptomatic unvaccinated class $I.$ 
	Individuals in the isolated compartment $Q$ are assumed not to contribute to disease transmission, as they are separated from the susceptible population.

	Apart from the assumption of  homogeneous mixing of the population the model is designed based on the following main assumptions:  
	\begin{enumerate}[(1)]
		\item  Individuals who recover from COVID--19 acquire temporary immunity and eventually return to the susceptible class after the immunity wanes. 
		\item   The severity of disease and COVID-19-induced mortality differ based on vaccination status. Specifically, \( \delta \) denotes the disease-induced death rate for unvaccinated individuals, while \( \delta_1 \) represents the corresponding rate for vaccinated individuals.
		\item Disease transmission is primarily driven by interactions between infectious and susceptible individuals. Consequently, the force of infection \( \lambda \) is formulated in terms of the infectious compartments and their respective contact rates with susceptibles.
		\item A constant recruitment rate \( \Lambda \) (due to births and migration) is assumed into the susceptible class \( S \). During lockdown periods, \( \Lambda \) is primarily attributable to natural births.
	\end{enumerate}
	
	Based on  these assumptions outlined above,  we now formulate the system of differential  equations governing the dynamics of each  compartment.
	In particular,  equations \eqref{S} and \eqref{V} describe  the  evolution of susceptible unvaccinated ($S$) and vaccinated ($V$) populations, respectively:%
	\begin{eqnarray} \label{S}
		S'(t)&=& \Lambda - [\lambda+(\sigma+\mu)]S +\varphi R,\\ 
		\label{V}
		V'(t)&=& \sigma S-[(1-\rho)\lambda  +\mu ]V+(\omega-\varphi)R.
	\end{eqnarray}
	Here, $\Lambda$ denotes the recruitment rate into the   $S$ class, which may result from births and migration. The parameter $\sigma$ represents rate at which susceptible individuals are vaccinated, while $\mu$ is the natural death rate. The vaccine effectiveness is denoted by    $\rho$, indicating the degree to which vaccination reduces susceptibility to infection. The parameter $\omega$ captures the overall rate at which individuals exit the $R$ class, and   $\varphi$ accounts for the fraction of recovered individuals who lose their temporary immunity and either refuse vaccination or lack access to it. 
	
	The dynamics of infectious classes   $A,\,I,\,A_1$ and $I_1$ are described by the system of  equations \eqref{eq A}--\eqref{eq I1}: 
		\begin{align}
			A'(t)&= \eta\lambda S  -(\mu+\theta+\gamma_1)A, \label{eq A}\\
			I'(t) &= (1-\eta)\lambda S  -(\epsilon+\gamma_2+\delta+\mu)I, \label{eq I}\\
			A_1'(t)&= \phi(1-\rho)\lambda V -(\theta_1+\gamma_4+\mu)A_1, \label{eq A1}\\
			I_1'(t) &=  (1-\phi)(1-\rho)\lambda V-(\epsilon_1+\gamma_5+\delta_1+\mu)I_1, \label{eq I1}
		\end{align}
	In this system, $\eta$ denotes  the proportion of unvaccinated individuals who become asymptomatic upon infection, while $\phi$ represents the corresponding proportion among vaccinated individuals. The parameters
	$\theta$ and $\gamma_1$ are the isolation/quarantine and recovery rates, respectively, for the unvaccinated asymptomatic class $A;$ and  $\theta_1$ and $\gamma_4$ are the their counterparts for the vaccinated asymptomatic class $A_1$.
	For the symptomatic classes, $\epsilon$ and $\epsilon_1$ denote isolation rates, $\gamma_2$ and $\gamma_5$ are recovery rates, and $\delta$ and $\delta_1$ are COVID-19-induced mortality rates for the unvaccinated $I$ and vaccinated $I_1$ compartments, respectively.
	
	The dynamics of the quarantine/isolation compartment $Q$ and recovered compartment $R$ are governed by     equations: 
		\begin{align}
			\label{Q}
			Q'(t) &= \theta A +\epsilon I+\theta_1 A_1+\epsilon_1 I_1 - (\gamma_3+\delta+\mu)Q, \\
			\label{R}
			R'(t) &= \gamma_1 A + \gamma_2I+\gamma_3 Q+\gamma_4 A_1+\gamma_5I_1 - (\omega + \mu)R.
		\end{align} 
	Here, $\gamma_3$ represents the  recovery rate from the isolation compartment $Q$ to the recovered compartment $R$.
	The complete transmission dynamics of the disease in the population  is illustrated in Figure~\ref{Figure1} and a summary description of the model's state variables and parameters   is   presented in Table~\ref{Table1}.
	
	Collecting equations \eqref{S}–\eqref{R}, we arrive at the following system of nonlinear ordinary differential equations governing the model dynamics:
	\begin{subequations} \label{eq:(2.1)}
		\begin{align}
			S'(t) &= \Lambda- [\lambda+\sigma+\mu]S +\varphi R, \label{eq:(2.1)S} \\
			V'(t) &= \sigma S-[(1-\rho)\lambda  +\mu ]V+(\omega-\varphi)R,\label{eq:(2.1)V} \\
			A'(t) &= \eta\lambda S  -(\theta+\gamma_1+\mu)A,\label{eq:(2.1)A}\\
			I'(t)&= (1-\eta)\lambda S -(\epsilon+\delta+\gamma_2+\mu)I\label{eq:(2.1)I},\\
			A_1'(t) &= \phi(1-\rho)\lambda V -(\theta_1+\gamma_4+\mu)A_1,\label{eq:(2.1)A1}\\
			I_1'(t) &= (1-\phi)(1-\rho)\lambda V-(\epsilon_1+\gamma_5+\delta_1+\mu)I_1,\label{eq:(2.1)I1}\\
			Q'(t)&= \theta A+\theta_1A_1 +\epsilon I+\epsilon_1I_1 - (\gamma_3+\delta+\mu)Q,\label{eq:(2.1)Q}\\
			R'(t) &= \gamma_1 A + \gamma_2I+\gamma_3 Q+\gamma_4A_1+\gamma_5I_1 - (\omega + \mu)R.\label{eq:(2.1)R}
		\end{align}
	\end{subequations}
	The system \eqref{eq:(2.1)} is equipped with   non-negative initial conditions: 
	\begin{eqnarray}\label{(2.4)}
		S(0)=S_0,\;V(0)=V_0,\;A(0) = A_0,\; I(0)=I_0,\; Q(0) = Q_0,\; \mbox{and} \;\;R(0)=R_0.
	\end{eqnarray}%
	
	 The next section is devoted to the analysis of model \eqref{eq:(2.1)} together with the initial conditions in \eqref{(2.4)}.
	
	\begin{figure}[ht!]
		\centering
		\begin{tikzpicture}[->,>=stealth',shorten >=1pt,node distance=2.5cm,auto,main node/.style={rectangle,rounded corners,draw,align=center}]
			\fill(-2,0) node[draw, rounded corners] (1) {$ S $}
			(-2,-3)node[draw, rounded corners] (2) {$ V $}
			(2,2) node[draw, rounded corners] (3) {$ A $}
			(1,1)node[draw, rounded corners] (4) {$ I $}
			(1,-4)node[draw, rounded corners] (5) {$ I_1 $} 
			(2,-5)node[draw, rounded corners] (6) {$ A_1 $}
			(4,-1.5) node[draw, rounded corners] (7) {$ Q $}
			(7,-1.5) node[draw, rounded corners] (8) {$ R $};
			\path
			(-3,0) edge node {{\footnotesize $\Lambda$}} (1);
			\path
			(1) edge node [swap] {\footnotesize $\sigma S$} (2) edge node [sloped, anchor=center, above] {\footnotesize $\eta\lambda S$} (3) edge node [sloped, anchor=center, below] {\footnotesize $(1-\eta)\lambda S$} (4) edge node [sloped, anchor=center, above] {\footnotesize $\mu S$} (-2.7,-0.7);
			\path
			(2) edge node [sloped, anchor=center, above] {\footnotesize  $(1-\phi)(1-\rho)\lambda V$} (5) edge node [sloped, anchor=center, below] {\footnotesize $\phi(1-\rho)\lambda V$} (6)edge node [sloped, anchor=center, above] {\footnotesize  $\mu V$}(-3,-3);
			\path
			(3) edge node [sloped, anchor=center, above] {\footnotesize  $\theta A$} (7) edge node [sloped, anchor=center, above] {\footnotesize $\gamma_1 A$} (8) edge node [sloped, anchor=center, above] {\footnotesize  $\mu A$}(1,2);
			\path
			(4) edge node [sloped, anchor=center, above] {\footnotesize $\epsilon I$} (7) edge node [sloped, anchor=center, above] {\footnotesize $\gamma_2I$} (8)  edge node [sloped, anchor=center, above] {\footnotesize  $(\mu+\delta)I$}(1,-1);
			\path
			(5) edge node [sloped, anchor=center, above] {\footnotesize $\epsilon_1I_1$} (7) edge node [sloped, anchor=center, above] {\footnotesize $\gamma_5I_1$} (8)edge node [sloped, anchor=center, below] {\footnotesize   $(\mu+\delta_1)I_1$}(1,-2);
			\path
			(6) edge node [sloped, anchor=center, above] {\footnotesize $\theta_1A_1$} (7) edge node [sloped, anchor=center, above] {\footnotesize $\gamma_4A_1$} (8) edge node [sloped, anchor=center, below] {\footnotesize  $\mu A_1$}(1,-5);
			\path
			(7) edge node [sloped, anchor=center, above] {\footnotesize $\gamma_3Q $} (8)edge node [sloped, anchor=center, above] {\footnotesize  $(\mu+\delta) Q$\qquad\qquad}(2,-1.5);
			\path
			(8) edge node [sloped, anchor=center, above] {\footnotesize  $\mu R$}(8,-1.5);
			\draw[->](8)--(7,3)--node[sloped, anchor=center, above]{{\footnotesize $\varphi R$}}(-2,3)--(1);
			\draw[->](8)--(7,-6)--node[sloped, anchor=center, below]{{\footnotesize $(\omega-\varphi) R$}}(-2,-6)--(2); 
		\end{tikzpicture}
		\caption{The schematic diagram for model \eqref{eq:(2.1)}.}
		\label{Figure1}
	\end{figure}

	\begin{table}  
		\begin{center}
			\begin{tabular}{cp{8cm}}
				\hline
				Variable&Description\\
				\hline
				$S$& Susceptible class \\
				$V$& Vaccinated class\\
				$A$ & Asymptomatic   and not vaccinated class\\
				$A_1$ & Asymptomatic   but vaccinated class\\
				$I$ & symptomatic   and not vaccinated class \\
				$I_1$ & Symptomatic  but vaccinated class\\
				$Q$& A class containing hospitalized individuals\\
				$R$& Recovered class\\
				\hline
				Parameter & Description \\ \hline
				$\Lambda $& Recruitment rate for $S$ class \\
				$\sigma$&  Rate of vaccination\\
				$\mu$& Natural death rate\\
				$\theta$& Rate of quarantine from class $A$  \\ 
				$\theta_1$& Rate of quarantine  from class $A_1$ \\
				$\gamma_1 ,\;\gamma_2, \;\gamma_3,\;\gamma_4,\;\gamma_5$& Rates of transfer   from  $A,\;I,\;Q,\;A_1$ and $I_1$ classes, respectively to $R$ class\\
				$\omega$& Rate of losing immunity from $R$ class \\ 
				$\varphi$& Fraction of recovered individuals who lost immunity but not vaccinated\\ 
				$\rho$& Vaccine effectiveness\\
				$\eta$ & Fraction of infected  individuals who remain asymptomatic  \\
				$\phi$& Fraction of infected but vaccinated individuals who remain asymptomatic\\
				$\epsilon$& Rate of isolation from $I$ class   \\
				$\epsilon_1$& Rate of isolation from $I_1$ class \\
				$\delta$ & Death rate due to COVID-19 for $I$ and $Q$ classes\\
				$\delta_1$ & Death rate due to COVID-19 for $I_1$ class\\
				$\beta$ & Effective contact rate in the community  \\
				$\nu,\;\nu_1, \;\kappa$ & Modification  parameters\\
				\hline
			\end{tabular}
			\caption{Description of  variables and parameters of model \eqref{eq:(2.1)}.} \label{Table1}
		\end{center}
	\end{table}

	\section{Quantitative and qualitative analysis} \label{analysis}

	 A critical first step in analyzing model \eqref{eq:(2.1)} is demonstrating its epidemiological feasibility, such as ensuring solutions remain biologically meaningful (e.g., non-negative and bounded). We formalize this property in Theorem \ref{Theorem3.1}, with the full proof provided in  \ref{Appendix A}. 
	
	\begin{theorem}  \label{Theorem3.1}
		The  model \eqref{eq:(2.1)} is a dynamical system on the  region  
		{
			\footnotesize 
				\begin{eqnarray}
					\Omega=\Big\{(S,V, A,I,A_1,I_1, Q,  R) \in \mathbb{R}_{+}^{8}:
					0\leq S+V+A+A_1+I+I_1+Q+R=N\leq \frac{\Lambda}{\mu}\Big\}, \label{eq maincone}
				\end{eqnarray}
		}  
		where $\mathbb{R}_+^8$ denotes the portion of $\mathbb{R}^8$ with all its components    non-negative.   
	\end{theorem}
	 Having established the feasibility of the system \eqref{eq:(2.1)}–\eqref{(2.4)}, we proceed to derive its equilibrium points and analyze their stability.
		The equilibria are obtained by setting the right-hand side of \eqref{eq:(2.1)} to zero. In particular, by setting all infected compartments to zero, we obtain the trivial equilibrium, 
			denoted by  $E_0,$ which represents  the disease-free equilibrium (DFE) given by \eqref{DFE}:
	\begin{eqnarray}\label{DFE}
		E_0=\left(S,V,A,I, A_1,I_1,Q,R\right)=\left(\frac{\Lambda}{\sigma+\mu}, \frac{\sigma\Lambda}{\mu(\sigma+\mu)},0,0,0,0,0,0\right) .
	\end{eqnarray} 
	If a non-trivial equilibrium solution $E^*$ exists, where at least one infected compartment is non-zero, it is referred to as an endemic equilibrium. We denote this equilibrium as  
	\begin{align}
		E^* = (S^*,\,V^*,\,A^*,\,I^*,\,A_1^*,\,I_1^*,Q^*,R^*) \label{eq:EE}
	\end{align}
	representing persistence infection within the population, even when the system is at equilibrium. We analyze the existence and stability of $E^*$ for the system \eqref{eq:(2.1)} in Subsection \ref{subsec: backward}.

  To analyze the stability of the equilibrium points, we compute the threshold parameter  known as the basic reproduction number, denoted by $\mathcal{R}_0$. Using the Next Generation Matrix (NGM) method as in \cite{Castilo_2002, Diekmann_2000, van_2008, Shuai_2013} we obtain
	\begin{align} \label{BRN}
		\mathcal{R}_0 = \mathcal{R}_A + \mathcal{R}_I + \mathcal{R}_{A_1} + \mathcal{R}_{I_1},
	\end{align}
	where 
	\begin{align*}
		\displaystyle \mathcal{R}_A = \frac{\nu B_1 }{k_1},\,~  \mathcal{R}_I = \frac{B_2}{k_2},\,~ \mathcal{R}_{A_1} = \frac{\nu_1 B_3}{k_3},\,~ \mathcal{R}_{I_1} = \frac{\kappa B_4}{k_4},
	\end{align*}
	and $B_i$ and $k_i$ for $i=1,\,2,\,3,\,4$ are defined as
	\begin{align*}
		&B_1=\frac{\eta \beta \mu}{\sigma+\mu},\;B_2=\frac{(1-\eta) \beta \mu}{\sigma+\mu},\;B_3=\frac{\phi(1-\rho) \beta \sigma}{\sigma+\mu},\;  \;B_4=\frac{(1-\phi)(1-\rho) \beta \sigma}{\sigma+\mu},\\
		&k_1=\theta +\gamma_1+\mu,\;k_2=\epsilon+\gamma_2+\delta+\mu,\;k_3=\theta_1+\gamma_4+\mu,\; 
		 k_4=\epsilon_1+\gamma_5+\delta_1+\mu . 
	\end{align*}
	  A detailed derivation of  $\mathcal{R}_0$ is provided in  \ref{BRN calculation}.

	\begin{remark}
		In \eqref{BRN}, the terms $\displaystyle \mathcal{R}_A  ,\,  \mathcal{R}_I,  \, \mathcal{R}_{A_1} $ and $    \mathcal{R}_{I_1} $ represent the  respective   contributions of the four infectious compartments $A,\, I,\, A_1$ and $I_1$  to the overall transmission potential of the diseases.
	\end{remark}
	
	\begin{remark}
		  When both vaccination and recovery confer permanent immunity (i.e., $\rho=1$ and $\omega = 0$), the basic reproduction number in \eqref{BRN} simplifies to 
		\begin{eqnarray}\label{BRNWV}
			\mathcal{R}_0=  \mathcal{R}_A+\mathcal{R}_I,
		\end{eqnarray}
		indicating that only the unvaccinated infectious compartments contribute to disease transmission under this assumption.
	\end{remark}
	The local asymptotic stability of the disease free equilibrium $E_0$, as stated in Theorem~\ref{Theorem3.3},  follows directly from    Theorem~2 of  \cite{van_2008}. 
	\begin{theorem} \label{Theorem3.3}
		The   disease-free equilibrium, $E_0$ of the  model \eqref{eq:(2.1)} is locally asymptotically stable   if  $\mathcal{R}_0 <1$  and unstable whenever $\mathcal{R}_0 > 1$.
	\end{theorem}

	\subsection{Existence of backward bifurcation}\label{subsec: backward}

	The  epidemiological implication of Theorem \ref{Theorem3.3}  is that  when $\mathcal{R}_0<1$, a small  influx of infected individuals into the community  will not trigger an  outbreak, and the disease eventually dies out. 
	  However, for disease elimination to be guaranteed regardless  of the initial  conditions, it is necessary to show that the   disease-free equilibrium ($E_0$) is globally asymptotically stable when $\mathcal{R}_0<1$. 
	  If the model exhibits a backward bifurcation,  then the disease may persist in the population even   when $\mathcal{R}_0<1,$  thereby  undermining this threshold behavior associated with the basic reproduction number. 
	
	To investigate the possibility of a backward bifurcation in the model, 
	 we examine the existence of a non-trivial endemic equilibrium
	 \(E^*= \left(S^*,V^*, A^*,I^*,A_1^*,I^*_1,Q^*,R^*\right)\) 
	 under the condition $\mathcal{R}_0<1$.  
	This equilibrium satisfies the  system of equations \eqref{eq:(2.3)}:
	\begin{equation} \label{eq:(2.3)} 
		\begin{aligned}
			&\Lambda- [\lambda^*+\sigma+\mu]S^* +\varphi R^*=0, \\&
			\sigma S^*-[(1-\rho)\lambda^*  +\mu ]V^*+(\omega-\varphi)R^*=0, \\ &
			\eta\lambda^* S^*  -k_1A^*=0,\\&
			(1-\eta)\lambda^* S^*-k_2I^*=0,\\&
			\phi(1-\rho)\lambda^* V^* -k_3A_1^*=0,\\&
			(1-\phi)(1-\rho)\lambda^* V^*-k_4 I_1^*=0,\\&
			\theta A^*+\theta_1A_1^* +\epsilon I^*+\epsilon_1I_1^* - k_5Q^*=0,\\&
			\gamma_1 A^* + \gamma_2I^*+\gamma_3 Q^*+\gamma_4A_1^*+\gamma_5I_1^* - (\omega + \mu)R^*=0,
		\end{aligned}
	\end{equation}
	where  \begin{eqnarray}\label{FFEE}
		\lambda^*=\beta \frac{I^*+\nu A^*+\nu_1A_1^*+\kappa I_1^*}{N^*}
	\end{eqnarray}
	is the force of infection and 
	\begin{eqnarray}\label{TPEE}
		N^*=S^*+V^*+A^*+I^*+A^*_1+I^*_1+Q^*+R^*,
	\end{eqnarray}
	the total population at the equilibrium.

	  While solving the full system \eqref{eq:(2.3)} is analytically challenging,  each equilibrium compartment, namely $A^*,\,I^*,\,A_1^*,\,I_1^*,\,Q^*,\,R^*,\, S^*$ and $V^*$, can be expressed  explicitly in terms of $\lambda^*$ as follows: 
	\begin{align}
		&A^*=\frac{\eta \lambda^* S^*}{k_1},\;\quad A^*_1=\frac{\phi (1-\rho)\lambda^* V^*}{k_3},\;\quad I^*_1=\frac{(1-\phi) (1-\rho)\lambda^* V^*}{k_4},\nonumber\\
		&Q^*=\frac{\lambda^*}{k_5}\left(t_1S^*+t_2 V^*\right),\;\quad ~R^*=\frac{\lambda^*}{\omega+\mu}\left(t_3S^*+t_4 V^*\right),\;\quad I^*=\frac{(1-\eta) \lambda^* S^*}{k_2}\nonumber \\
		& S^*=t_5 V^*+t_6,\,\quad \mbox{and}\quad V^*=\frac{(\sigma+x t_3 \lambda^*)t_6}{(1-\rho)\lambda^*+\mu-(\sigma t_5+x\lambda^*(t_3t_5+t_4))} \label{25}
	\end{align} 
	where %
	\begin{align*}
		&  k_5 = \gamma_3+\delta+\mu, \quad t_1= \frac{\eta \theta}{k_1}+\frac{\epsilon(1-\eta)}{k_2},\quad t_2= \frac{\phi(1-\rho) \theta_1}{k_3}+\frac{\epsilon_1(1-\phi)(1-\rho)}{k_4}, \quad \\
		&t_3=\frac{\gamma_1 \eta}{k_1}+\frac{\eta_2(1-\eta)}{k_2}+\frac{\gamma_3 t_1}{k_5}, \quad t_4=\frac{\gamma_3 t_2}{k_5}+\frac{\gamma_4\phi(1-\rho)}{k_3}+\frac{\gamma_5(1-\phi)(1-\rho) }{k_4},\quad  \\
		&  t_5=\frac{\varphi t_4}{(\omega+\mu)(\lambda^*+\sigma+\mu)-\varphi\lambda^* t_3},\quad
		t_6= \frac{(\omega+\mu)\Lambda}{(\omega+\mu)(\lambda^*+\sigma+\mu)-\varphi\lambda^* t_3} ,   \quad  x= \frac{\omega-\varphi}{\omega+\mu}.
	\end{align*}

	  Through algebraic manipulation of equations (\ref{FFEE}), (\ref{TPEE}), and (\ref{25}), we obtain a polynomial of degree four in $\lambda^*$  of the form:%
	\begin{eqnarray}\label{244}
		H(\lambda^*)=\lambda^{*}P(\lambda^*),
	\end{eqnarray}
	where
	\begin{eqnarray}\label{26}
		P(\lambda^*)= \mathcal{Q}_3 (\lambda^*)^3+ \mathcal{Q}_2(\lambda^*)^2+\mathcal{Q}_1\lambda^*+\mathcal{Q}_0,
	\end{eqnarray}
	and
	\begin{eqnarray}
		\mathcal{Q}_3&=&F_3(\omega+\mu-\varphi t_3)>0,\nonumber\\
		\mathcal{Q}_2&=&F_3(\omega+\mu)(\sigma+\mu)+F_2(\omega+\mu-\varphi t_3),\nonumber\\
		\mathcal{Q}_1&=&F_2(\omega+\mu)(\sigma+\mu)+F_1(\omega+\mu-\varphi t_3),\nonumber\\
		\mathcal{Q}_0&=&F_1(\omega+\mu)(\sigma+\mu)\,=\,(\omega+\mu)(\sigma+\mu)^2(1-\mathcal{R}_0). \label{eq:polycoefficients}
	\end{eqnarray}
	The auxiliary terms $F_1,\,F_2,\,F_3$ are defined as:
	\begin{eqnarray*}
		F_1&=&\sigma D_1+\mu D_3,\\
		F_2&=&D_1t_3x+\sigma D_2+\left(1-\rho-t_4 x\right)D_3+\mu D_4,\\
		F_3&=&D_2 t_3 x+ \left(1-\rho-t_4 x\right)D_4,
	\end{eqnarray*}
 with the intermediate coefficients $D_1,\,D_2,\,D_3,\,D_4$ given by:
	\begin{eqnarray*}
		D_1&=&1-\left(\frac{\nu_1\phi(1-\rho)\beta}{k_3}+\frac{\kappa (1-\phi)(1-\rho)\beta}{k_4}\right),\\
		D_2&=&\frac{\phi(1-\rho)}{k_3}+\frac{(1-\phi)(1-\rho)}{k_4}+\frac{t_2}{k_5}+\frac{t_4}{\omega+\mu},\\
		D_3&=&1-\left(\frac{\nu \eta \beta}{k_1}+\frac{(1-\eta)\beta}{k_2}\right),\\
		D_4&=& \frac{ \eta }{k_1}+\frac{(1-\eta)}{k_2}+\frac{t_1}{k_5}+\frac{t_3}{\omega+\mu}.
	\end{eqnarray*}
	 
	 The equation $H(\lambda^*)$ in \eqref{244} admits the trivial solution $\lambda^*=0$, which corresponds to the disease-free equilibrium ($E_0$) defined in \eqref{DFE}. Any additional non-negative real roots of  $P(\lambda^*)=0$, when they exist, correspond to endemic equilibrium points of the system.
	From the coefficient expressions in \eqref{eq:polycoefficients}, we observe that $\mathcal{Q}_0>0$ when $\mathcal{R}_0<1$, and $\mathcal{Q}_3$ remains positive for all parameter values.
	
	To investigate the possibility of multiple endemic equilibria, we apply  Descartes' rule of signs to the cubic polynomial $P(\lambda^*)$ given in \eqref{26}. 
	This classical method enables  us to  systematically determine the number of positive real roots  by examining sign changes in the sequence of polynomial coefficients under various parameter regimes. 
	In particular, the signs of $Q_1,\,Q_2$ and the value of $\mathcal{R}_0$ play a critical role in determining the root structure of  $P(\lambda^*)$.
	Table \ref{Table2} summarizes the different cases, classifying the number of possible positive real roots based on these key parameters. 
	\begin{table}[ht!]
		\centering
		\begin{tabular}{p{1cm}p{0.5cm}p{0.5cm}p{0.5cm}p{0.5cm}p{1cm}p{1.5cm}p{4cm}}\hline
			Cases &$\mathcal{Q}_3$&$\mathcal{Q}_2$&$\mathcal{Q}_1$&$\mathcal{Q}_0$&$\mathcal{R}_0$& No of sign changes & No of possible equilibrium (roots)\\
			\hline
			1&+&+&+&+&$<1$&0&0\\
			&+&+&+&-&$>1$&1&1\\\hline
			2&+&+&-&+&$<1$&2&0, 2\\
			&+&+&-&-&$>1$&1&1\\\hline
			3&+&-&+&+&$<1$&2&0, 2\\
			&+&-&+&-&$>1$&3&1, 3\\\hline
			4&+&-&-&+&$<1$&2&0, 2\\
			&+&-&-&-&$>1$&1&1\\\hline
		\end{tabular}
		\caption{Number of possible positive roots for (\ref{26}).} \label{Table2}
	\end{table} 

	  Theorem \ref{Theorem3.33} follows directly from the classification of root structures summarized in Table~\ref{Table2}. 
	
\begin{theorem}\label{Theorem3.33}
		
	Consider the model \eqref{eq:(2.1)}. Then,
	\begin{enumerate}[(1)]
		\item the model admits a unique endemic equilibrium if cases 1, 2, or 4 occur with $\mathcal{R}_0>1.$
		\item  the model may admit multiple endemic equilibria if case 3 occurs with $\mathcal{R}_0>1.$
		\item  the model may admit multiple endemic equilibria  if cases 2, 3, or 4 occur with $\mathcal{R}_0<1.$
		\item  the model admits no endemic equilibrium if case 1 occurs with $\mathcal{R}_0<1.$
	\end{enumerate}

	\end{theorem}
	
	\begin{remark}
		Item (3) of Theorem~\ref{Theorem3.33}  establishes the co-existence of disease-free equilibrium and one or more endemic equilibria when $\mathcal{R}_0<1$. This implies that the system \eqref{eq:(2.1)}  may exhibit a backward bifurcation phenomenon, where the disease can persist in  the population even though the basic reproduction number is below the classical threshold of unity.  
		In such cases, reducing  $\mathcal{R}_0$ below one is not sufficient to guarantee disease elimination, and additional control measures may be required to steer the system towards the disease-free-equilibrium.
	\end{remark}

\begin{figure}[h!] 
		\centering
			\includegraphics[scale=0.5]{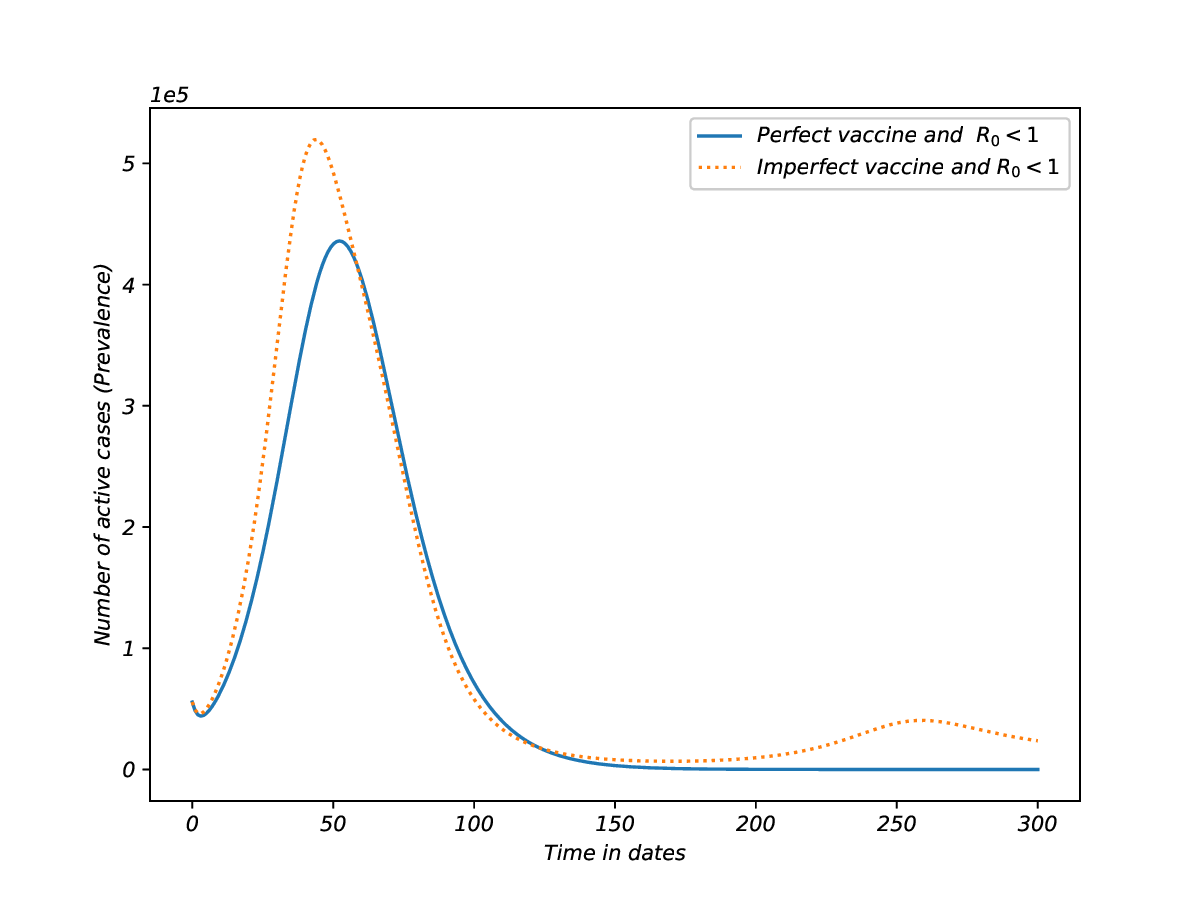} 
%
		\caption{Stability analysis of the disease-free equilibrium (DFE) comparing perfect vaccination ($\rho=1,\,\omega=0$)with imperfect vaccination. Although both scenarios satisfy $\mathcal{R}_0,$ the system under perfect vaccination (blue curve) shows global DFE stability, while imperfect vaccination (orange curve) leads to bistability—coexistence of DFE and endemic equilibria—indicating backward bifurcation. This illustrates how vaccine efficacy ($\rho$) and waning immunity ($\omega$) critically affect eradication thresholds.}
			\label{fig:sigma0rho1}
\end{figure}

	 The phenomenon of a backward bifurcation in the model (\ref{eq:(2.1)}) is   formally established in Theorem \ref{Theorem3.77}.
		\begin{theorem} \label{Theorem3.77}
			The model (\ref{eq:(2.1)}) exhibits a backward bifurcation at  $\mathcal{R}_0=1$. %
		\end{theorem}
		The proof of Theorem \ref{Theorem3.77}, which relies on the center manifold theory and  Theorem 4.1 of \cite{Castilo_2004},   is provided in  \ref{Appendix BB}. Furthermore, the existence of a backward bifurcation  for $\mathcal{R}_0 < 1 $ is demonstrated numerically   in Figure~\ref{fig:sigma0rho1}.

	When there is  no reinfection (i.e. $\omega=0$) and the vaccine is perfect (i.e.  $\rho=1$),  the model  (\ref{eq:(2.1)}) does not exhibit a backward bifurcation.
	Under these conditions, the disease-free equilibrium is globally asymptotically stable for   $\mathcal{R}_0<1$,  as stated in Theorem \ref{Theorem3.4}, with the proof   provided in \ref{Appendix B}.
	
	\begin{theorem} 
		Under the assumption $\omega=0$ and $\rho=1$, the disease-free equilibrium, $E_0$ of the  model \eqref{eq:(2.1)} is globally asymptotically stable  whenever $ \mathcal{R}_0 <1$.
		\label{Theorem3.4}
	\end{theorem}
	Numerical simulation illustrating the   global asymptotic stability of the disease-free equilibrium for $\mathcal{R}_0 < 1 $ under the conditions $\omega=0$ and $\rho=1$ is presented in  Figure~\ref{fig:sigma0rho1}. 

	  When reinfection is absent ($\omega=0$) and vaccines confer perfect immunity ($\rho=1$), the coefficients in the polynomial equation \eqref{26} simplify such that  $\mathcal{Q}_3=0,\;\mathcal{Q}_2>0$ and $\mathcal{Q}_1>0$.  According to Table~\ref{Table2}, this configuration ensures a unique positive root of \eqref{26}. Consequently, the system \eqref{eq:(2.1)} admits a unique endemic equilibrium $E^*;$
	\begin{eqnarray}\label{30}
		E^{*}=(S^{*},V^{*},A^{*},I^{*},0,0,\,Q^{*}, R^{*})
	\end{eqnarray}
	for $ {\mathcal{R}_0}>1$.  
	This implies that under these conditions,   the model does not exhibit a backward bifurcation, even when the basic reproduction number exceeds unity. In this case, the total population size $N$,  the force of infection   $\lambda$ and the basic reproduction $\mathcal{R}_0$ simplify to 
	\begin{eqnarray}
		{N}&=& S+V+A+I+Q+R,\\
		{\lambda}&=&\beta\frac{I+\nu A}{ {N}},\\
		{\mathcal{R}_0}&=& \frac{\nu B_1}{k_1}+\frac{B_2}{k_2}.
	\end{eqnarray}

	Consequently,   the following remark follows from the preceding discussion and Theorem 4.1, item (iv), of \cite{Castilo_2004}.
	%

		\begin{remark}\label{Theorem3.7}
	For the model \eqref{eq:(2.1)}, a unique endemic equilibrium exists and is locally asymptotically stable, provided that $\mathcal{R}_0 > 1$ and $\mathcal{R}_0$ is sufficiently close to unity.
	\end{remark}

	 The next theorem establishes that perfect vaccination (i.e, $\rho=1 $) and the absence of reinfection (i.e, $\omega=0$) can suppress disease transmission, even when   $\mathcal{R}_0>1.$

\begin{theorem}\label{Theorem3.8}
	For $\rho=1$ and $\omega=0,$ let $\;I^* = I^*(\sigma)$ denote the infectious component of the endemic equilibrium corresponding to the vaccination rate parameter $\sigma$. Then  we have 
	$$I^*(\sigma) < I^*(0),$$
	that is, increasing vaccination rate reduces the infectious population at equilibrium.
\end{theorem}

\begin{proof}
	When $\rho=1$, the endemic equilibrium point is given by
	 $$E^*=(S^*,V^*,A^*,I^*,0,0,Q^*,R^*),$$ 
	 and the conservation law  (equation \eqref{Gron}) becomes
	\[\Lambda-\mu N^*-\delta(I^*+Q^*)=0,\] 
	where $N^*=S^*+V^*+A^*+I^*+Q^*+R^*$. 
	Rearranging, we get
	\[(\mu+\delta)I^*=\Lambda-\mu(S^*+V^*+A^*+R^*)-(\mu+\delta)Q^*.\] 
	Thus, the infectious component at equilibrium as a function of $\sigma$ is
	\[ I^* = I^*(\sigma) = \frac{\Lambda-\mu(S^*+V^*+A^*+R^*)-(\mu+\delta)Q^*}{\mu+\delta}.\]
	Since the vaccination rate $\sigma$ affects $V^*$ but not other components, we have
	\begin{eqnarray*}
		I^*(\sigma)&=&\dfrac{\Lambda-\mu(S^*+V^*+A^*+R^*)-(\mu+\delta)Q^*}{\mu+\delta} \\
		&\leq&\dfrac{\Lambda-\mu(S^*+A^*+R^*)-(\mu+\delta)Q^*}{\mu+\delta}=I^*(0).
	\end{eqnarray*}
	This completes the proof. 
\end{proof}

\section{Model Fitting} \label{fitting}

This study develops an epidemiological model to assess COVID-19 vaccine effectiveness, calibrated and validated using case data from South Africa. Although the analysis focuses on South Africa, the vaccines deployed in the county, such as Pfizer-BioNTech, Johnson \& Johnson, Oxford-AstraZeneca, among others, were also central to global vaccination efforts, spanning high-income countries (e.g., United States, United Kingdom) to middle-income nations (e.g., Brazil, India).

Despite delayed rollouts in Africa and other low-income regions, many of which relied heavily on  WHO's COVAX initiative for vaccine supply, the same vaccines were deployed as in wealthier countries. 
However, real-world effectiveness varied substantially across populations due to factors like variant prevalence, healthcare infrastructure, and demographic differences \cite{rubin2021covid,andrews2022covid,chung2021effectiveness,chen2025dynamics}. 
Therefore, while our findings should not be generalized as universal vaccine performance metrics, the model remains applicable to countries with similar vaccination strategies, such as Botswana and Namibia. %

	COVID-19 case data for South African were obtained from the Johns Hopkins University COVID-19 repository \cite{JH_data}. The analysis follows a   three phase approach:
	\begin{enumerate}[(1)]
		\item Data Presentation, providing an overview of calibration datasets including case counts and vaccination timelines. 
		\item Parameter Estimation, with initial parametrization  based on WHO and CDC guidelines, as well as  published literature. 
		\item Model Fitting, where parameter estimates are  refined through least-squares optimization implemented via the Python \emph{lmfit} library.
	\end{enumerate}   %

\subsection{About the data}\label{sec: about the data}

\begin{figure}
	\centering
	\begin{subfigure}[b]{0.3\textwidth}
		\includegraphics[width=\textwidth]{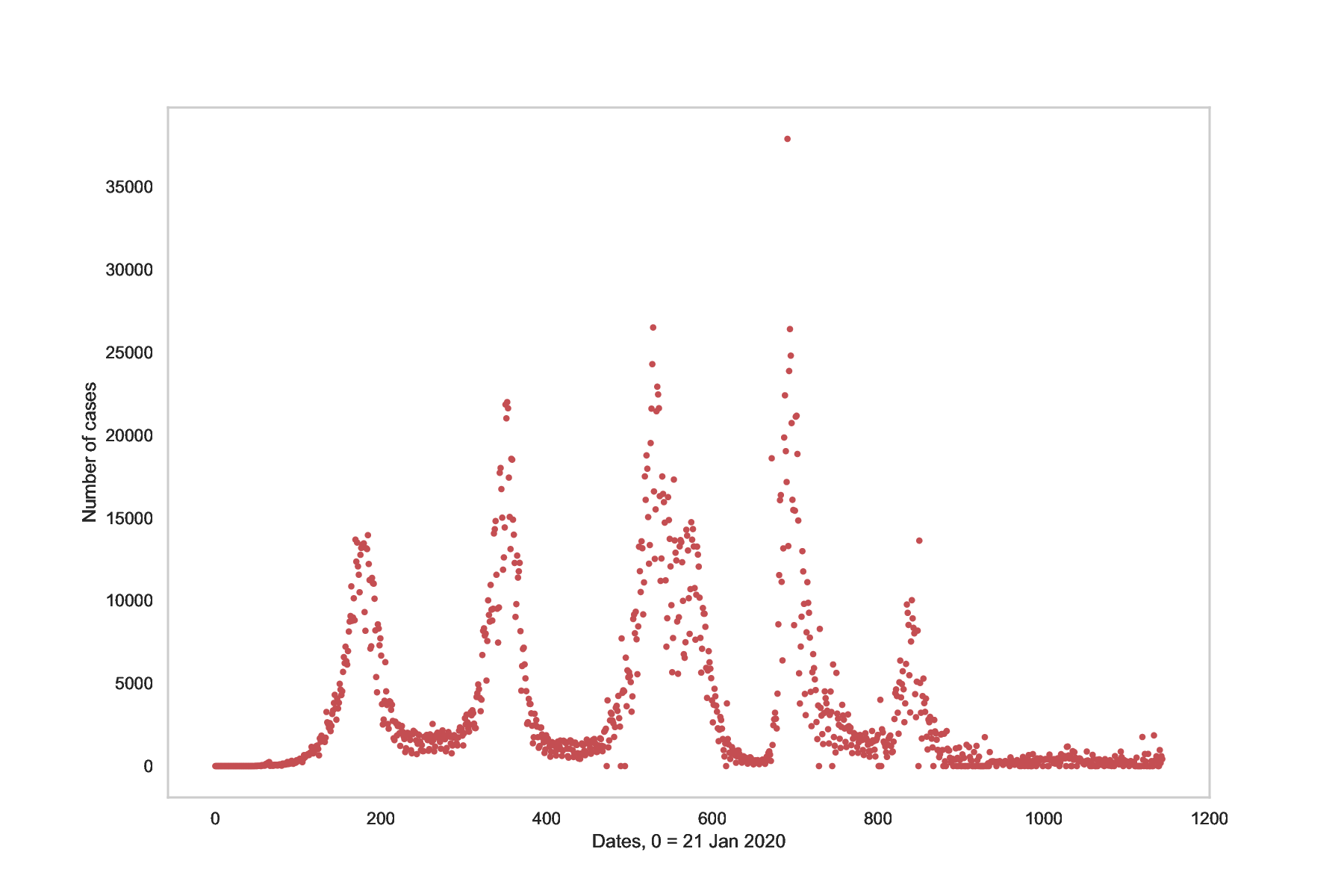}
		\subcaption{South African Daily infection data} \label{fig: dailyinfectiondotted}
	\end{subfigure}
	\begin{subfigure}[b]{0.3\textwidth}
		\includegraphics[width=\textwidth]{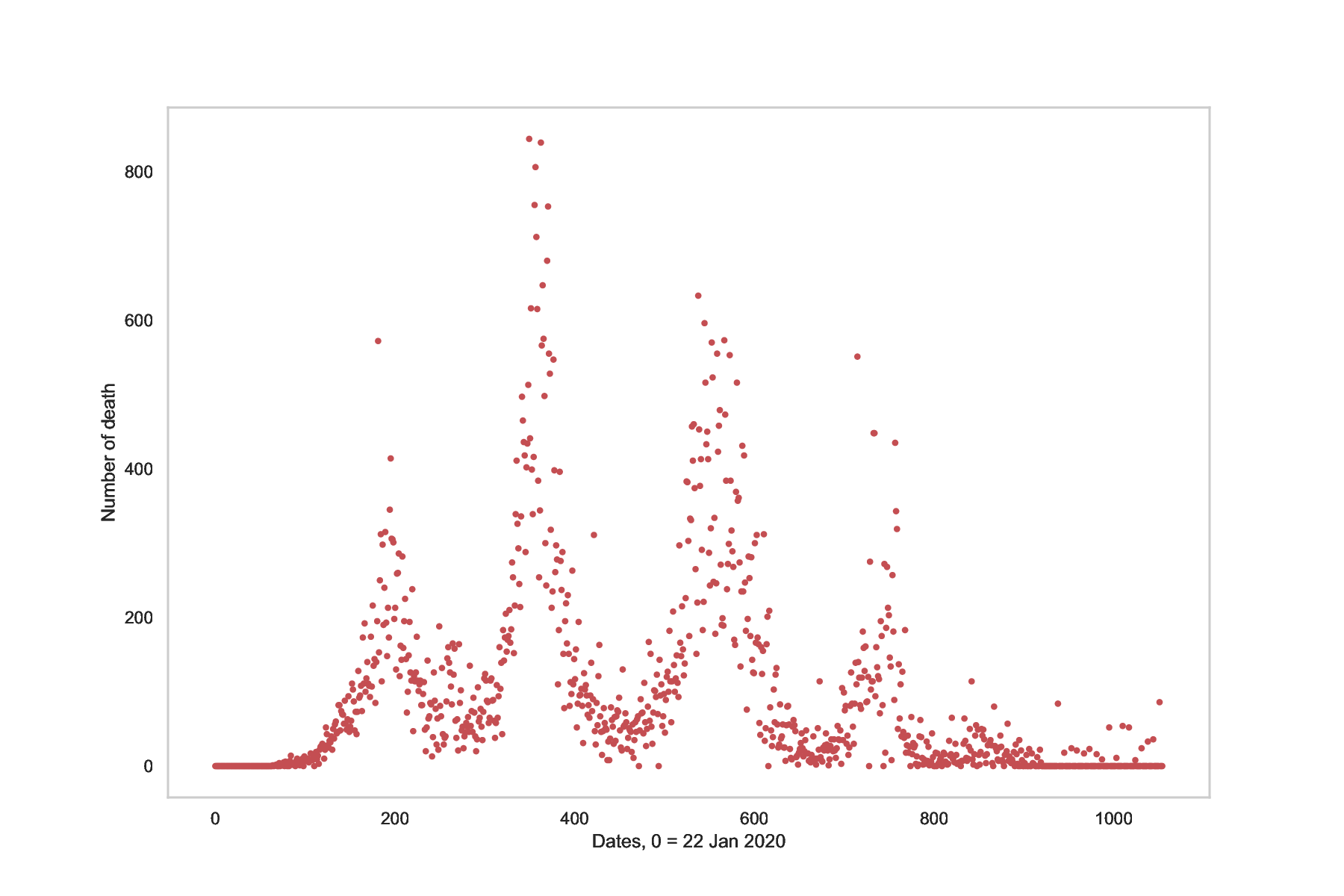}
		\subcaption{South African Daily death data} \label{fig: dailydeathdotted}
	\end{subfigure}
	\begin{subfigure}[b]{0.3\textwidth}
		\includegraphics[width=\textwidth]{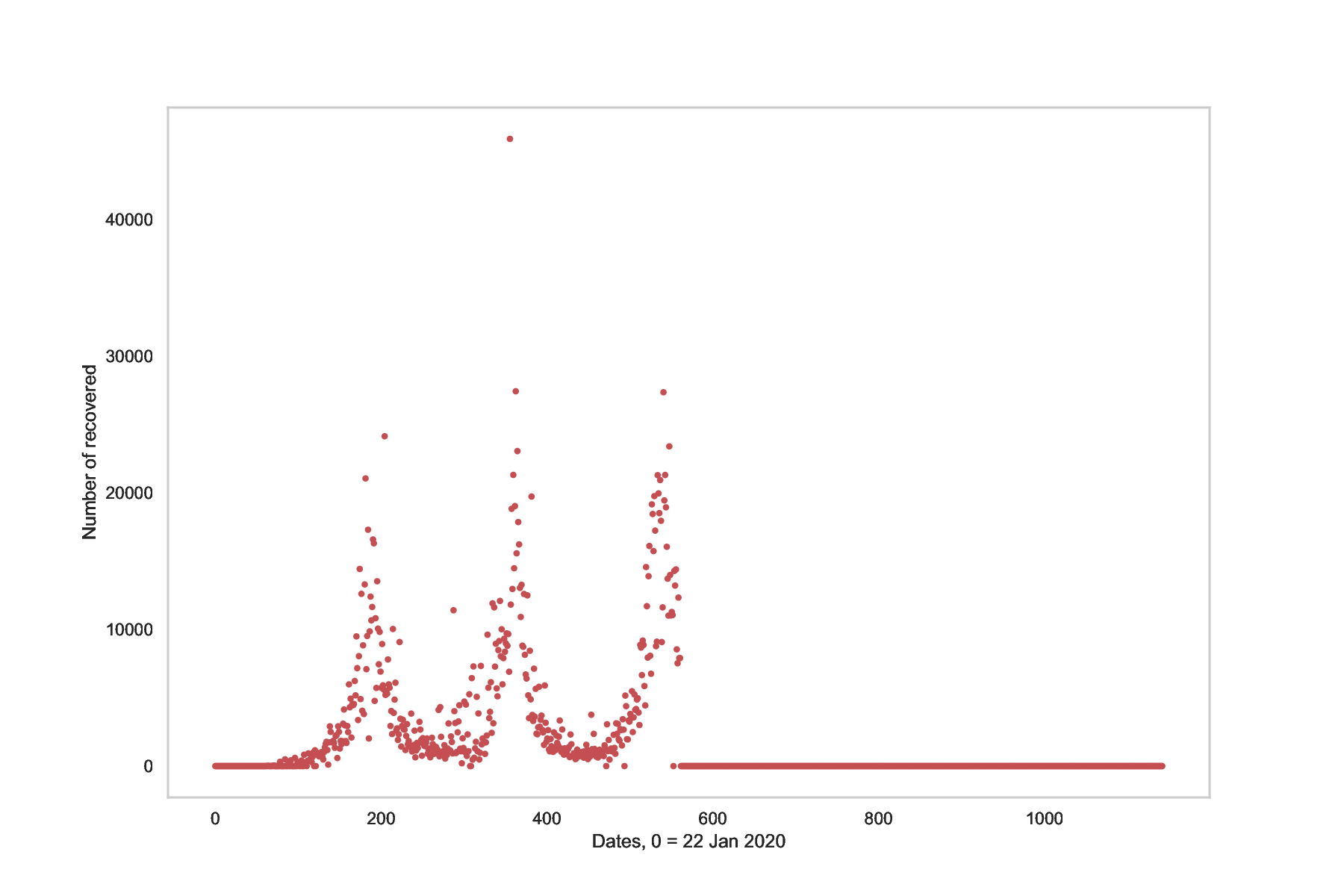}
		\subcaption{South African Daily recovery data}\label{fig: dailyrecoverydotted}
	\end{subfigure}
	\caption{South African daily COVID-19 data from \cite{JH_data}}
	\label{fig: SA daily covid data}
\end{figure}

\begin{figure}[h!]
	\centering
	
	\begin{subfigure}[b]{0.48\textwidth}
		\includegraphics[width=\textwidth]{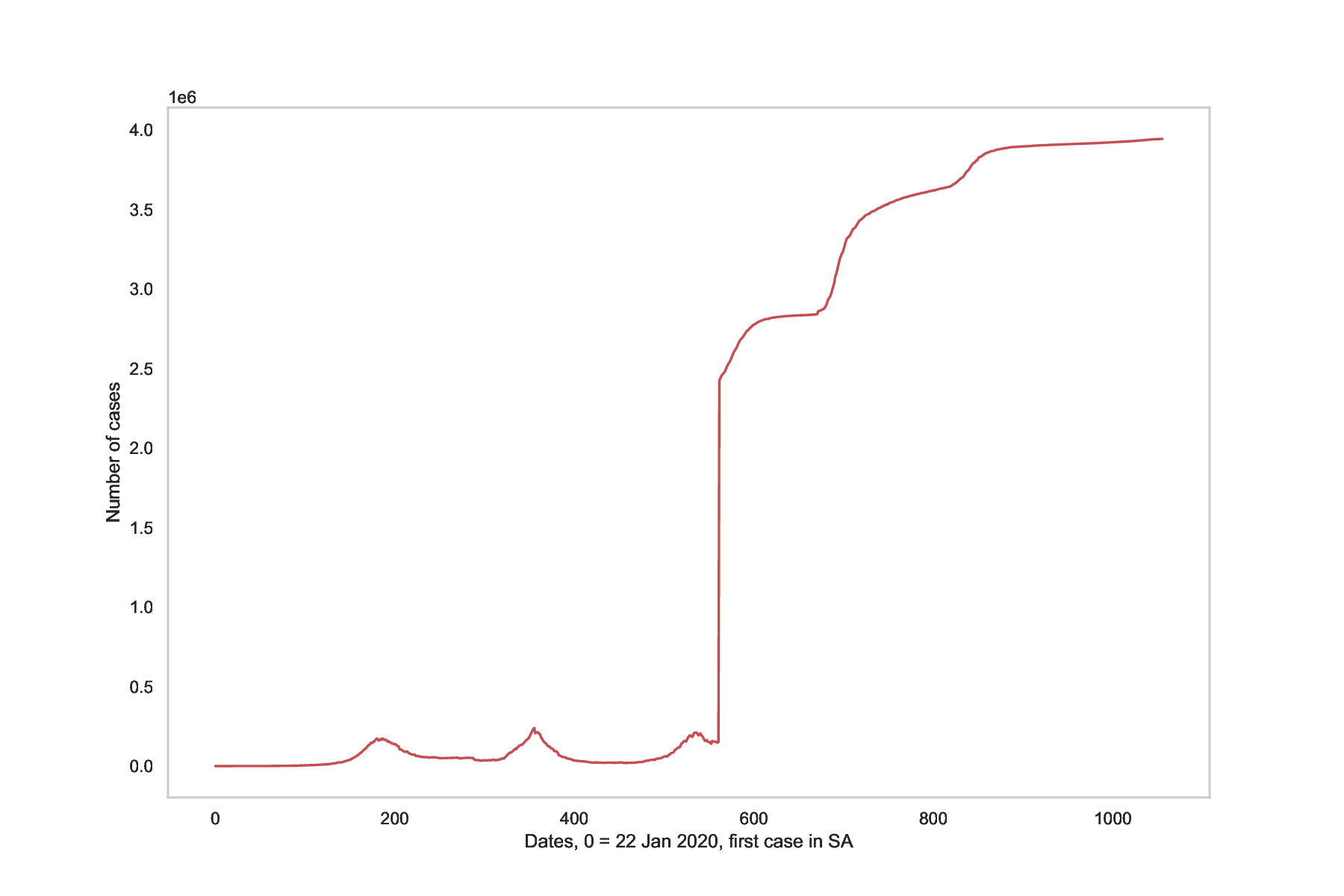}
		\subcaption{Active South African COVID19 cases from 22 January 2020 to 18 November 2022. } 
		\label{fig: active sa full}
	\end{subfigure}
	\begin{subfigure}[b]{0.48\textwidth}
			\includegraphics[width=\textwidth]{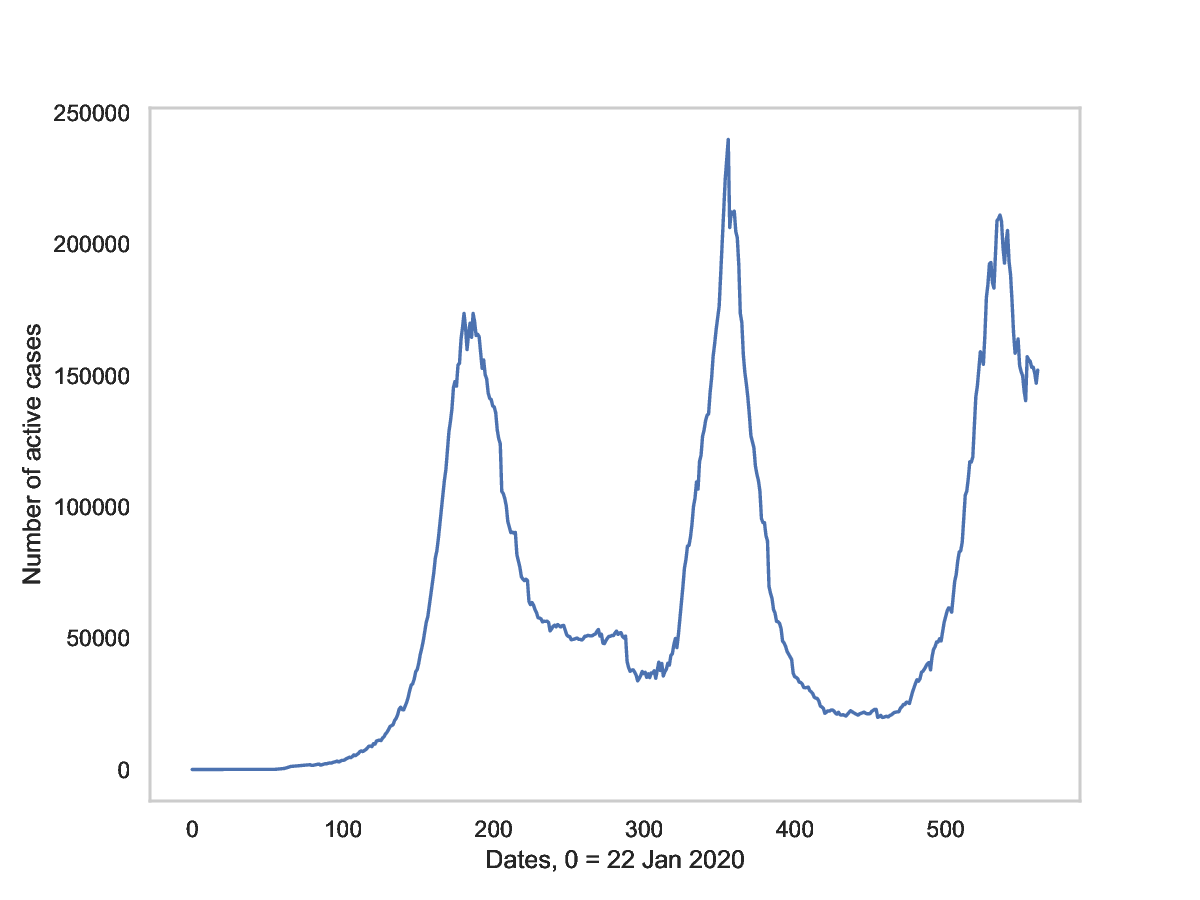} 
		\subcaption{ Active South African cases from 22 Jan 2020 to 5 August 2021.}
		\label{fig: sa active 561}
	\end{subfigure}
	\begin{subfigure}[b]{0.48\textwidth}
		\includegraphics[width=\textwidth]{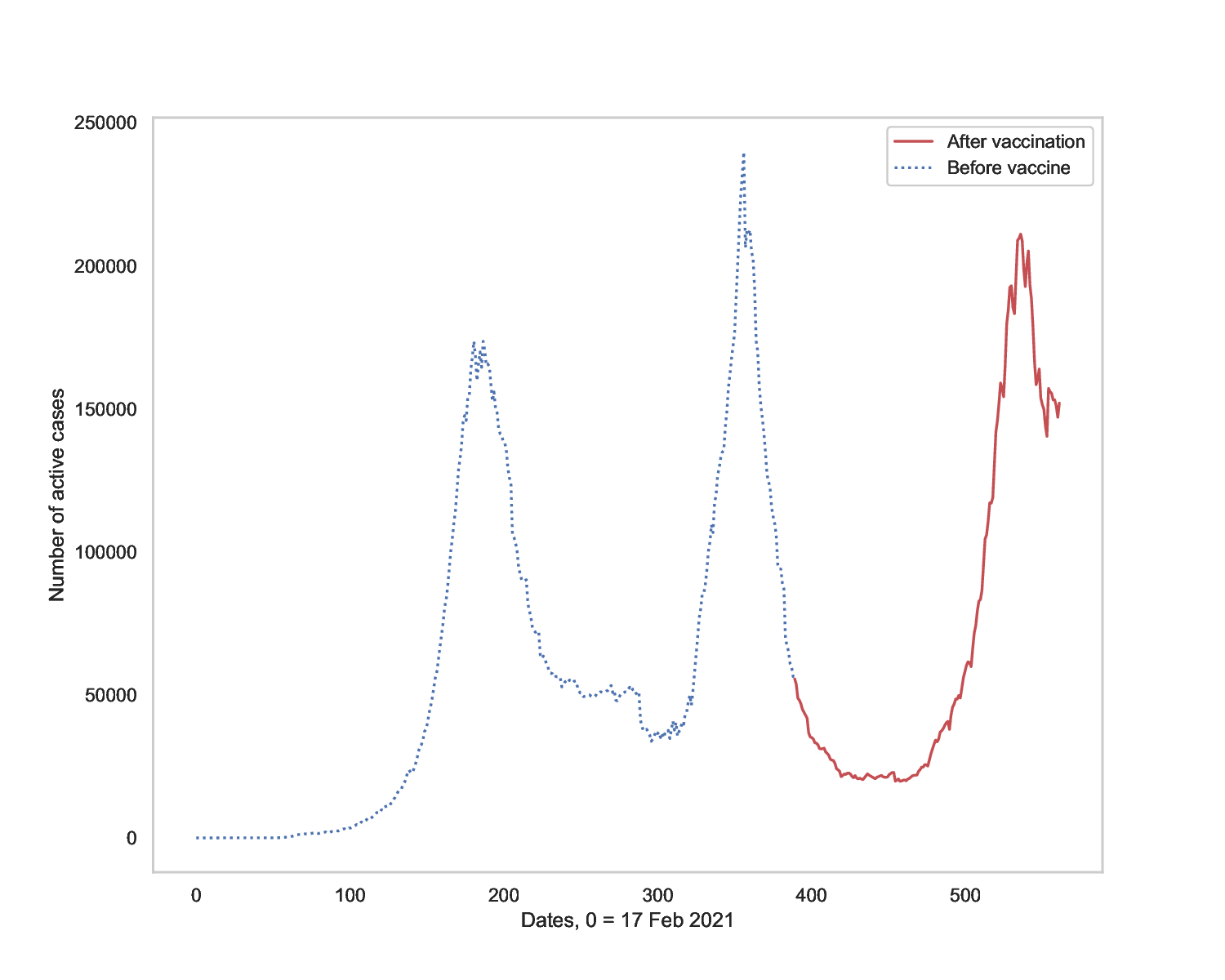} 
		\subcaption{The blue portion represents the data before vaccines introduced in South Africa }
		\label{fig: active blue-red}
	\end{subfigure}
	\begin{subfigure}[b]{0.48\textwidth}
		\includegraphics[width=\textwidth]{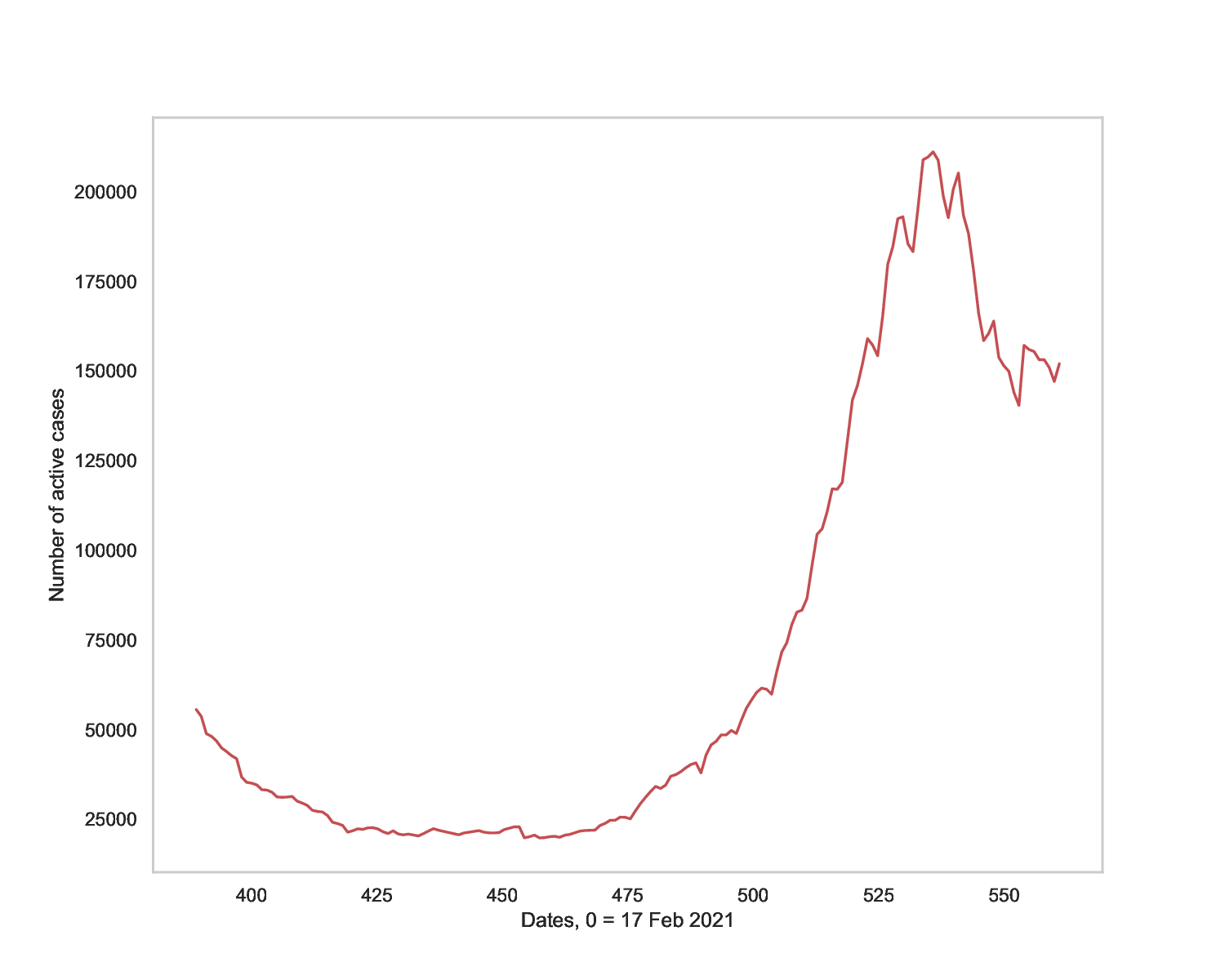} 
		\subcaption{Active South African cases between 17 Feb 2021 and 5 August 2021.}
		\label{fig: sa active 172}
	\end{subfigure}
	\caption{South African COVID-19 data}
	\label{fig: prevalence data}
\end{figure}

Global COVID-19 data, publicly available from Johns Hopkins University \cite{JH_data}, were downloaded and processed using Python.
 Figure \ref{fig: SA daily covid data}   presents South Africa's  daily COVID-19 statistics, including reported daily infections, recoveries, and deaths. 
 During data processing,    we identified a discontinuity in the daily recovery records beginning on day 561 (relative to January 22, 2020, the date of South Africa’s first reported case), which corresponds to  August 5, 2021. 
 This anomaly results in an artificial spike in active cases on day   561, as illustrated in Figure \ref{fig: active sa full}. 
 To ensure consistency in subsequent analysis, the dataset was truncated at     day 561 (see Figure \ref{fig: sa active 561}). %

\begin{figure}
	\begin{subfigure}[b]{0.48\textwidth}
		\includegraphics[width=\textwidth]{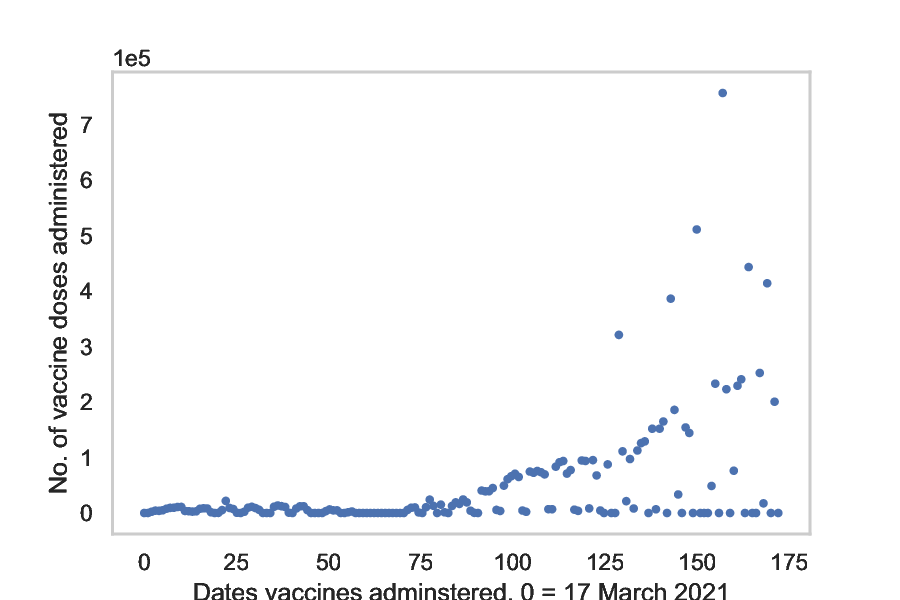}
		\subcaption{Daily vaccination time series}
		\label{fig: sa vaccine}
	\end{subfigure}
	\begin{subfigure}[b]{0.48\textwidth}
		\includegraphics[width=\textwidth]{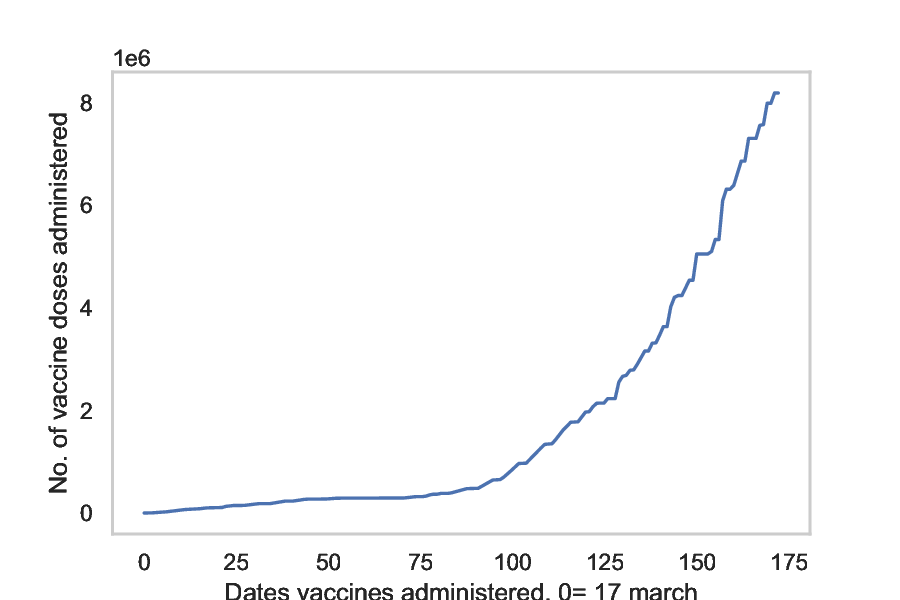} 
		\subcaption{Total vaccination time series}
		\label{fig: sa vaccine daily total}
	\end{subfigure}
	\caption{South African vaccination data for the first 172 days.}
	\label{fig: sa vaccine data}
\end{figure}

The South African vaccine data were also obtained from the same Johns Hopkins University database \cite{JH_data}.
The first COVID-19 vaccines were administered in South Africa on February 17, 2021. 
Since our goal is to analyze the  impact of vaccinations, we focus on the data recorded between February 17, 2021 and August 5, 2021 as this period contains a complete record  of vaccination data and active cases data. This subset of the daily vaccination data is shown in Figure~\ref{fig: sa vaccine data}.
Figure~\ref{fig: active blue-red} highlights the periods  before and after the introduction of vaccines, while Figure~\ref{fig: sa active 172} presents the daily active cases for the first 172 days since the start of vaccination campaign (February 17 - August 5, 2021,).

\subsection{Initial Parameter Estimation} \label{initial parameter estimation}

To set our initial parameter estimations, we used data and literature from various credible sources.

According to the latest World Bank estimate \cite{{WorldBank}}, the South African population is approximately 60.2 million, with an annual growth rate of $1.2\%$ and annual death rate of $9.468/1000$. This yields  a daily recruitment rate $$\Lambda\approx3583.3$$ and daily non-COVID related death rate of $$\mu\approx 2.6433\times10^{-5}.$$ 
The daily vaccination rate, $\sigma$, is estimated by the daily average of vaccine doses administered per day between   February 17 and  August 5, 2021  (Figure \ref{fig: sa vaccine data}). During this 172 day period, South Africa administered 8182380 doses of vaccines doses. Thus, the estimated  vaccination rate is $$\sigma \approx 7.9\times10^{-4}.$$

It is well known that non of the COVID-19 vaccines  provided full protection against infection, but they significantly reduced symptom severity   and improved recovery outcomes \cite{CDCvaccine,selamawit}. 
Early reports considerable variation in vaccine effectiveness. According to \href{https://www.webmd.com/vaccines/covid-19-vaccine/covid-vaccines-compared}{WebMed} \cite{webofmedcine}, the efficacy of the WHO approved vaccines ranged from approximately 50\% to   95\%; see also \cite{sadoff2021safety,soheili2023efficacy,medicalnewstoday}.
We adopted the midpoint of this range, $\rho=0.75$ as our initial estimate for vaccine effectiveness. 
To define a plausible fitting range, we considered the possibility that real-world effectiveness may be lower than reported due to uncontrolled circumstances, such as variant prevalence, delayed dosing, or storage issues. Based on this, we selected the interval   $(0.3,1)$ as our fitting range for $\rho.$


To estimate the proportion  of symptomatic and asymptomatic infection in the population, we use a result from the case study conducted on a specific workplace setting  \cite{Paleker2021asymptomatic}. 
According to this study,   36.6\%  of the infected individuals remain asymptomatic through out their infection, while up to 45\% are  asymptomatic in the first few days.  
Based on these findings, we assume that 45\% of all infections are asymptomatic. Among these,  81.3\% remain asymptomatic throughout their infection, while the   remaining 18.7\% develop mild symptoms that prompt them to seek testing within one week, subsequently transitioning to   compartment $Q.$
Accordingly, we set the initial estimates as $\eta=0.45 $ and $ \theta =0.187/7 = 0.0267$, with fitting ranges of $(0.3,0.6)$ for $\eta$ and $(0.01,0.03)$   for  $\theta$. 

The authors were unable to identify a publicly available study that specifically estimates the asymptomatic-to-symptomatic infection proportion among vaccinated individuals. 
However, based on reports from  the World Health Organization and other sources,  vaccines are known to reduce the severity of symptoms and the risk of fatal outcomes. 
Consequently,  we assume that the proportion of asymptomatic infections among vaccinated individuals,  $\phi$, is greater than that of the general population, i.e.,  $\phi>\eta. $
In particular, we set 
$\phi \approx 0.5$  as an initial estimate and define the fitting range $(\eta, 1)$. 

According to \cite{Benisek2023} and the Department of health of South Africa, individuals who test positive for COVID-19 are required to isolate for a minimum of 10 days. 
Based on this guideline, we assume  that  individuals in compartment $Q$ remain in isolation for approximately   10--15 days. 
Considering that the overall   COVID-19 has a recovery rate is estimated at  97.4\%  \cite{worldometer}, we  define the fitting range for the recovery rate parameter $\gamma_3$ as $\displaystyle   (0.0694, 0.0974)$ and set an initial estimation of   $\displaystyle \gamma_3 = 0.09$.

According to \cite{bai2020presumed}, the recovery period for asymptomatic patients ranges from   1-2 weeks,  with the median duration of 9 days.
Based on this, we  assume a fitting range of  $\displaystyle   (0.0544 , 0.1167)$ for the recovery rate of asymptomatic individuals, $\gamma_1$ and set the initial estimate to $\gamma_1 = 1/9\times 0.813 = 0.0903.$

Furthermore, according to \cite{worldometer}, the overall recovery rate in South Africa  is about 97.4 \%, 
and  according to \cite{coviddeath},  the median time to death for non-survivors is  18.5 days.  Based on these reports, we define the fitting range for the death rate $\delta$ as $\displaystyle  (0.0012,0.0016)$ as our fitting range  and 
\(\delta = 0.0015,\)
as an initial estimate.

Given that vaccination reduces the risk of death, we assume that the death rate among vaccinated individuals, $\delta_1,$ is lower than that of the general population, i.e., $\delta_1<\delta$. Therefore, we set the fitting range for $\delta_1$ as $\displaystyle  (0,\delta)$   and take   $\delta_1 = 0.0011$ 
as an initial estimate. 

Next, we  estimate the parameters $\epsilon, \epsilon_1,\gamma_2, \gamma_4,\gamma_5, \text{ and }\theta_1.$ 

Hypothetically, every individual in compartment $I$ is expected to    experience some level of discomfort or pain. As a result,  most of them would either choose to self-isolate or use government facilities to protect their loved ones or require hospitalization due to the severity of symptoms. 
We assume that only 25\% of these individuals recover with out entering isolation, while the majority transition  to compartment $Q$.
Assuming that such transitions occur with in a week, we define the fitting range for $\epsilon$ as $\displaystyle   (0.1051, 0.1472)$  and set the initial estimate as $\epsilon = 0.1252$. 

According to \cite{walsh2020sars,zhou2020clinical}, the recovery period for symptomatic individuals ranges from   8 to 37 days, with a median of approximately 20 days. Based on this, we define the fitting range for the recovery rate parameter $\gamma_2$ as $\displaystyle   (0.0066,0.0313)$ and use  $\gamma_2=  0.0125$ as an initial estimate. 

To estimate the remaining parameters, such as $\theta_1, \epsilon_1, \gamma_4 $ and $ \gamma_5$, we  assume individuals in compartment    $I_1$ transition to compartments  $Q$ and $R$ in equal proportion. However, due to  vaccination, they are expected to recover quickly than their unvaccinated counterparts. 
Following a similar rationale to that used for compartment  $I$, and assuming that individuals self-isolate within 2 to 7 while those who do not isolate recover with in   8 to 38 days, we define the fitting ranges as   $\displaystyle\epsilon_1\in (0.1057 , 0.37)$ and $\displaystyle\gamma_5\in (0.0066 , 0.0313)$. We set the initial estimates as  
   $\epsilon_1 = 0.125 ,\, \gamma_5 = 0.0275$.

Finally, by employing an approach analogous to the one used for estimating   $\theta$ and $\gamma_1$, we determine the fitting ranges $\displaystyle\theta_1\in (0.01,0.1)$ and $\displaystyle\gamma_4\in (0.0544,0.1167)$, with initial estimates 
$\theta_1= 0.0267$ and $ \gamma_4 = 0.095 $.
A summary of all parameter estimates and their corresponding fitting ranges is  provided  in table \ref{tab: Table3}.

\subsection{Fitting the model to the South African   data collected on days 17 Feb 2021 - 5Aug 2021}  \label{sec: fitting}

\begin{figure}[h!]
	\centering
	\includegraphics[scale=0.4]{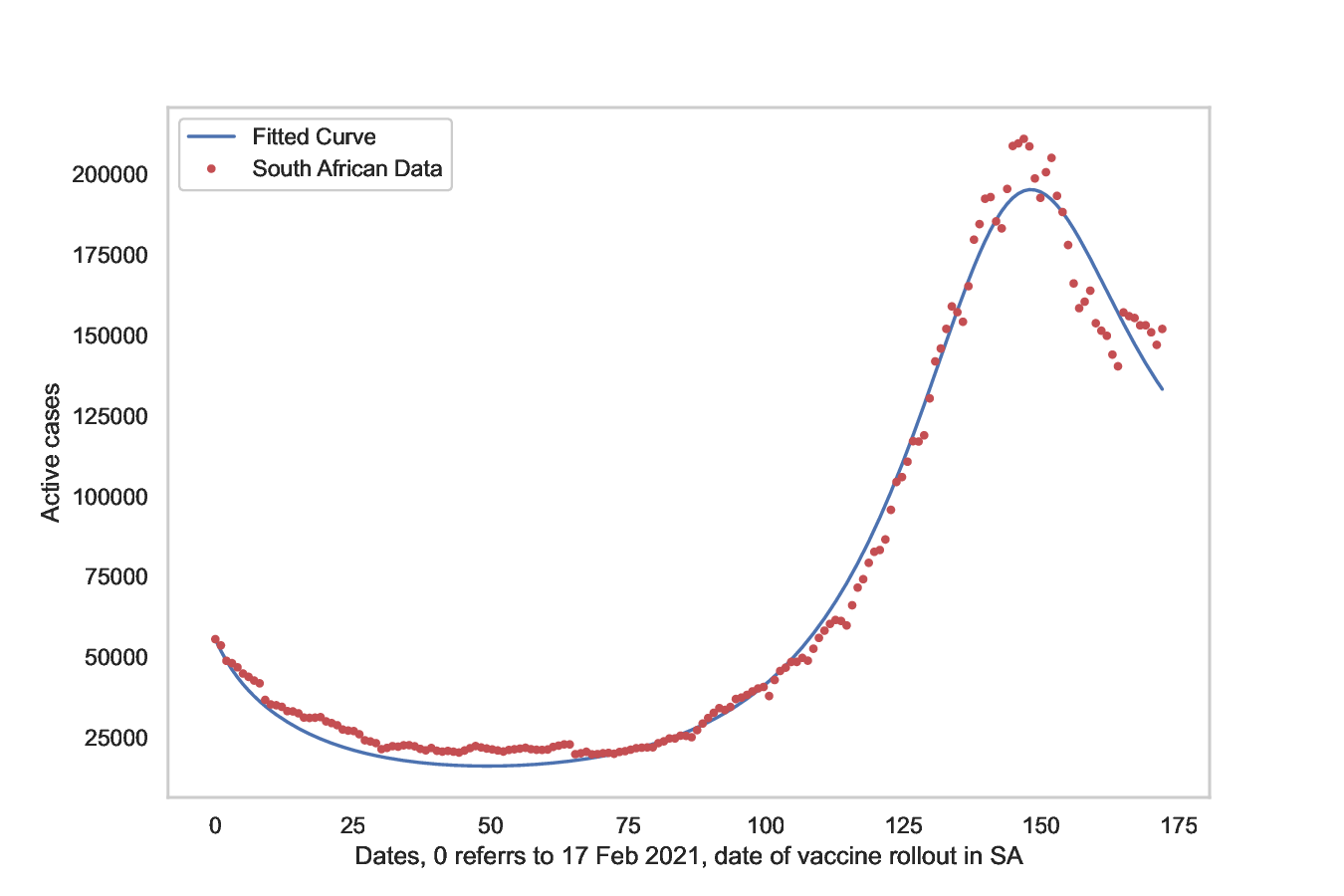}
	\caption{The model fitted to the South African Covid-19 Data from 17 Feb  to  5 Aug 2021.}
	\label{fig: fitted}
\end{figure}
Although our parameter estimation is primarily data-driven and supported by relevant literature, further refinement within the specified range is necessary through model fitting.
 The model is calibrated using   the data shown in Figure~\ref{fig: sa active 172}, which corresponds the relevant time period, i.e., February 17 - August 5, 2021. 
The fitting procedure is carried out using the open-source Python library \emph{lmfit}.
Table~\ref{tab: Table3} summarizes both  the initial parameter estimates, as outlined in Section \ref{initial parameter estimation}, and their optimized values obtained through the fitting process. 
The results of the fitting are shown in Figure \ref{fig: fitted}, which illustrates that the model effectively captures the observed transmission dynamics.

Since the primary goal of this analysis is to assess vaccine effectiveness through a mathematical model approach, particular emphasis is placed on parameters related to vaccination, namely, $\rho$, $\delta_1$, $\gamma_4$, and $\gamma_5$. 
The parameters $\gamma_4$ and $\gamma_5$ represent the rate at which vaccinated individuals transition to the recovered compartment $(R)$, while $\delta_1$ denotes the COVID-19-induced mortality rate among vaccinated individuals.
Lastly, the parameter $\rho$ quantifies the effectiveness of vaccines in preventing infection.

	According to  WHO, CDC, and other sources \cite{rubin2021covid,andrews2022covid,chung2021effectiveness}, vaccine effectiveness varies based on several factors, including the number of doses administered, the type of vaccine, the time elapsed since vaccination, the characteristics of the vaccinated population, and, most critically, the COVID-19 variant. 
	Reports from the WHO and the National Institute for Communicable Diseases of South Africa indicate that  the beta variant was dominant from February to May 2021, while the Delta variant was dominant from June to August 2021. For further details on variant prevalence in South Africa, see \cite{tegally2021detection}, and for insights into incorporating COVID-19 variants into epidemiological models and their impact on vaccine effectiveness, refer to \cite{chen2025dynamics} and the references in there.
These studies report vaccine efficacy ranging from as high as 95\% to as low as 50\%, depending on the variant and other influencing factors.
Importantly, such reported efficacy generally reflects the vaccines combined ability   to prevent symptomatic infection and to reduce the risk of hospitalization and death.

However, a recent report from the National Center for Immunization and Respiratory Diseases (NCIRD), published on February 1, 2024 \cite{NCIRD} highlights  improvements in COVID-19 vaccine performance.
Based on data collected from mid-September 2023 to January 2024, the report states that individuals vaccinated with recent vaccine formulations were about 54\% less likely to get infected with the virus. 
Given this, we anticipate that  our model fitting  will yield an estimated effectiveness rate for preventing infection  that is below 54\%.


In our model fitting, particular attention was given to the parameter  $\rho,$ which represents the   effectiveness of vaccines in preventing infection. 
Although, in theory,  $\rho$ can range from 0 to 1, it is unreasonable to assume that vaccines with 0\% effectiveness would  have received regulatory approval or widespread distribution.
Taken into account the urgency and uncertainty present during the initial stages of vaccine rollout,   we allowed $\rho$ to vary between 0.3 and 1, as specified in Subsection~\ref{initial parameter estimation}. This range accommodates potential variations in vaccine performance, even though existing  literature reports effectiveness rates typically ranging  from 50\% to 95\%. 

The fitted result suggests that $\rho \approx 0.3000$, which is lower than the reported effectiveness of currently produced vaccines \cite{NCIRD}, as anticipated.
This outcome aligns with the understanding that $\rho$ in our model specifically captures vaccine effectiveness in preventing infection, rather than the broader efficacy metrics reported in clinical studies.  
Additionally, the fitted parameters   suggest $\gamma_2 = 0.0066 $ and $\gamma_5 = 0.0232$, indicating that individuals in  compartment $I_1$ (vaccinated and infected)  recover more rapidly than those in compartment  $I$ (unvaccinated and infected). This supports the reports from WHO and CDC that vaccination reduces disease severity.
Furthermore, the results yield  $\delta = 0.0012$ and $\delta_1 = 2.02\times 10^{-14}\approx 0$, highlighting the significant protective effect of vaccination in reducing COVID-19 related mortality. These findings underscore the vital role of vaccines in saving lives during the pandemic.

Thus, these results support the conclusion that vaccines have played a significant role in reducing both disease severity and mortality.
This is evident from the relationship   $\gamma_2 \ll \gamma_5$, indicating faster recovery among vaccinated individuals,   from  $\delta_1 \lll \delta$, reflecting a dramatic reduction in COVID-19 induced mortality due to vaccination.  
However, since  $\rho < 0.5$, the effectiveness of vaccines in preventing infection appears to be lower than initially anticipated. But, as mentioned earlier, this could also be due to the new delta variant that appeared immediately after the introduction of vaccines \cite{chen2025dynamics}.
 For further discussion and   concluding remarks,  the reader is referred to Section~\ref{conclusion}.

\section{Local sensitivity analysis} \label{sensitivity}

In this section,  we analyze  the  sensitivity of the basic reproduction number  $\mathcal{R}_0$ to  slight changes in the parameters involved in equation \eqref{BRN}.  
This is achieved by computing  the partial derivative of $\mathcal{R}_0$ with respect to a parameter of interest, say $p$. 
The resulting quantity $$\frac{\partial \mathcal{R}_0}{\partial p}$$ is referred to as the local sensitivity index of $\mathcal{R}_0$ with respect  $p,$ and is denoted by $\displaystyle \gamma_p^{\mathcal{R}_0}$ (see  \cite{martcheva2015}), i.e., 
\begin{equation*}
	\gamma_p^{\mathcal{R}_0} =  \frac{\partial \mathcal{R}_0}{\partial p}.
\end{equation*}

for a better reflection of the correlation between $\mathcal{R}_0$ and parameter $p$, we   consider the normalized   sensitivity index,   denoted by $\displaystyle \epsilon^{\mathcal{R}_0}_p$, which is defined as 

\begin{equation*}
	\epsilon_{p}^{\mathcal{R}_0} := \frac{\partial \mathcal{R}_0}{\partial p} \frac{p}{\mathcal{R}_0}.  
\end{equation*} 

This normalized index expresses the relative changes in $\mathcal{R}_0$ resulting from relative changes in $p.$ 
In particular, $\epsilon_{p}^{\mathcal{R}_0} $ quantifies   the percentage  change in $\mathcal{R}_0$  
when $p$ is changed by a certain percent, say $ y\,\%$ (see \cite{martcheva2015,terefe2018}). 
That  means that
\begin{align*}
	\Delta\mathcal{R}_0 \% = \epsilon_{p}^{\mathcal{R}_0} y\%,
\end{align*}
or equivalently,
\begin{equation*}
	\epsilon_{p}^{\mathcal{R}_0}   = \frac{\Delta \mathcal{R}_0 \,\%}{\Delta p\,\%}. 
\end{equation*}


The local sensitivity index of $\mathcal{R}_0$ is summarized in Table~\ref{tab: Table3}, which shows the relationship between each of the parameter and $\mathcal{R}_0.$
 A visual representation of these indices is provided in  Figure~\ref{fig:sensgraph} for a visual presentation of the indices.
A negative sensitivity index, $\displaystyle \epsilon^{\mathcal{R}_0}_p$, indicates  that $\mathcal{R}_0$ and the parameter $p$ are inversely related, while a positive  sensitivity index implies a direct relationship. 
Moreover, the absolute value of $\epsilon_p^{\mathcal{R}_0}$ reflects the relative influence of parameter $p$ on the basic reproduction number  $\mathcal{R}_0$: larger values suggest greater sensitivity.

Since our goal is to assess the contribution of vaccination, we focused particularly  parameters associated with vaccination, such as $\rho,\,\phi,\, \theta_1,\, \epsilon_1,\, \gamma_4$ and $\gamma_5.$
Using the basic reproduction number $\mathcal{R}_0$ as a response function, several key conclusions can be drawn from the sensitivity indices summarized in  Table \ref{tab: Table3} and visualized Figure \ref{fig:sensgraph}.

Figure \ref{fig:sensgraph} highlights the   significance of parameters $\rho,\; \phi, \;\theta_1,\; \gamma_4,\; \epsilon_1, \; \beta,\; \nu_1$ and $\kappa. $
Notably, the relatively large sensitivity index $\epsilon^{\mathcal{R}_0}_\phi$ compared to $\epsilon^{\mathcal{R}_0}_\eta$  suggests  that vaccinated individuals with asymptomatic infections contributed substantially to ongoing transmission. 
In addition, the sensitivity indices  $\epsilon^{\mathcal{R}_0}_{\theta_1}$ and $\epsilon^{\mathcal{R}_0}_{\epsilon_1}$ underscore the importance of identifying and isolating infected individuals, even those who are vaccinated, in mitigating disease spread.

The parameter, $\rho$ which represents vaccine effectiveness of vaccines in preventing infection,  is inversely related to $\mathcal{R}_0$.
Thus, improving $\rho$ would significantly reduce disease transmission. 
The relatively small value of $\rho$, combined with  the large $|\epsilon^{\mathcal{R}_0}_{\rho}|$, underscores the importance of continued research aimed at developing more effective vaccines or enhancing the efficacy of existing ones. 

To further understand the role of vaccination, we compare some of the parameters with their unvaccinated counterparts. 
From Table~\ref{tab: Table3}, we observe that 
\begin{align*}
	|\epsilon^{\mathcal{R}_0}_{\theta}|<|\epsilon^{\mathcal{R}_0}_{\theta_1}| \qquad \text{and}\qquad |\epsilon^{\mathcal{R}_0}_{\epsilon}|<|\epsilon^{\mathcal{R}_0}_{\epsilon_1}|.
\end{align*} 
This suggests that identifying and isolating vaccinated individuals who become infected may have an even greater impact on reducing transmission than similar efforts among unvaccinated individuals. 
These findings highlight the need for continued monitoring of vaccinated individuals, along with public awareness campaigns that emphasize the imperfect nature of vaccine protection and encourage regular testing to determine infection status. 

\begin{figure}
	\begin{subfigure}[b]{0.48\textwidth}
		\includegraphics[width=\textwidth]{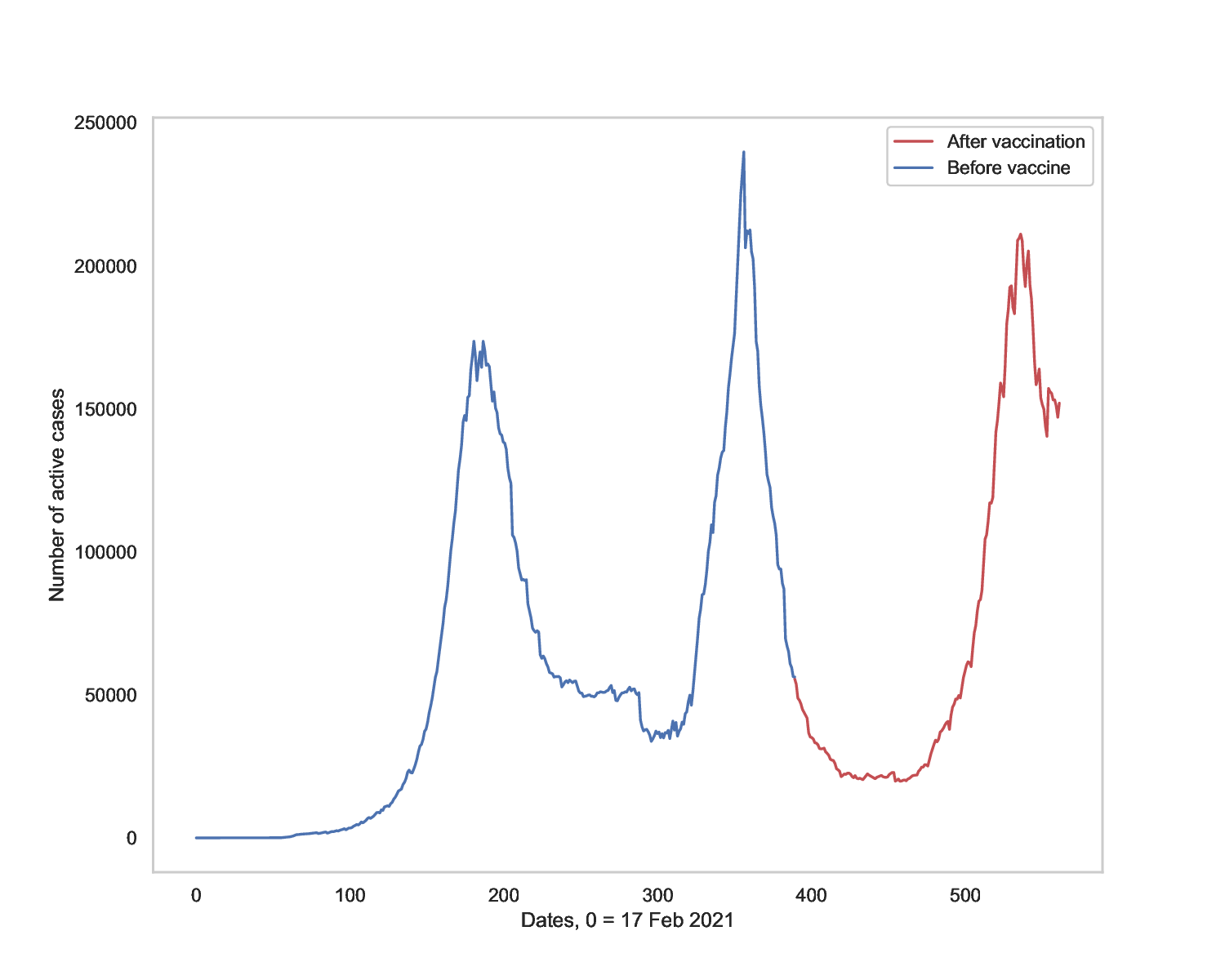}
		\subcaption{South African Active cases until 5 August 2021}
		\label{fig: sa active cases}
	\end{subfigure}
	\begin{subfigure}[b]{0.48\textwidth}
		\includegraphics[width=\textwidth]{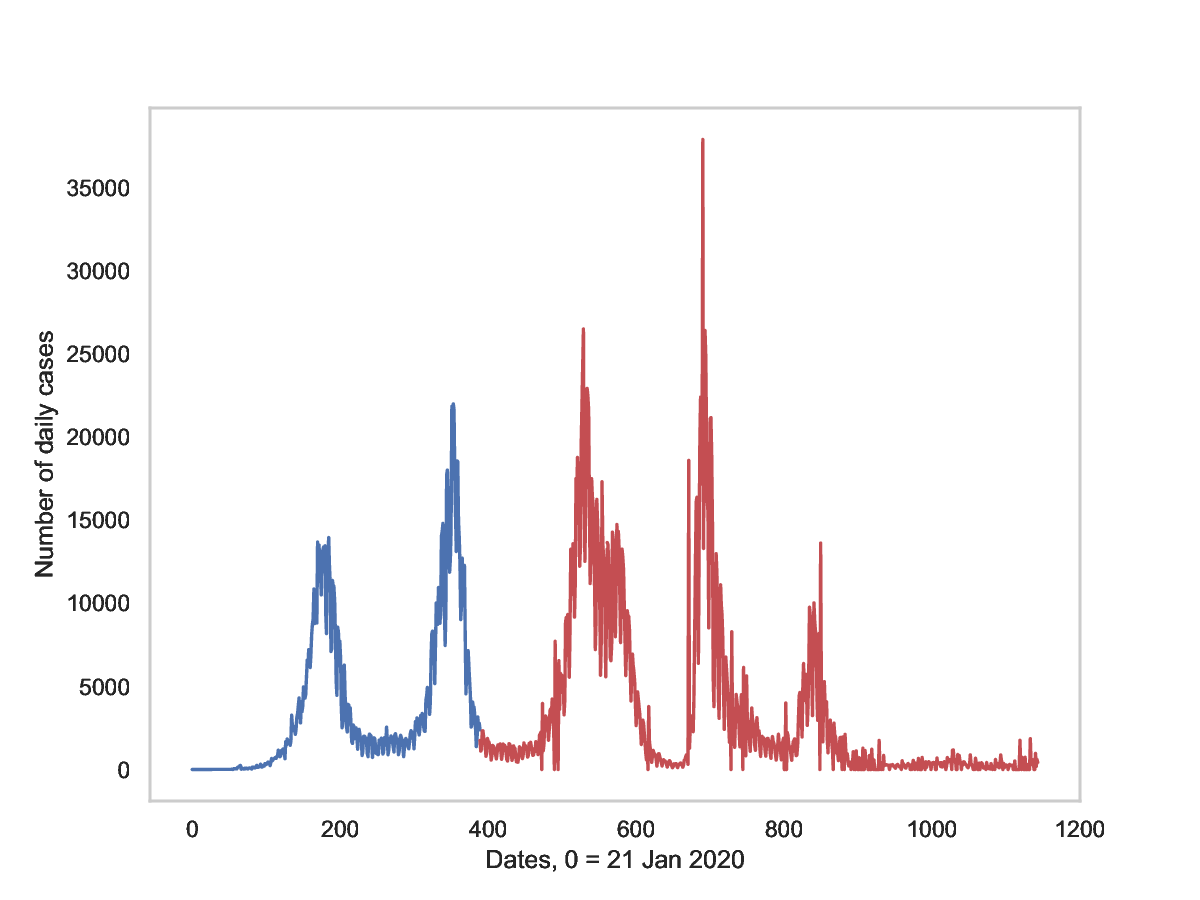} 
		\subcaption{South African daily infection cases until recently}
		\label{fig: sa ddaily infection}
	\end{subfigure}
	\caption{South African Active and daily infection data.}
	\label{fig: prev_vs_daily_infection}
\end{figure}

\begin{figure}[h!]
	\centering
	\includegraphics[scale=0.5]{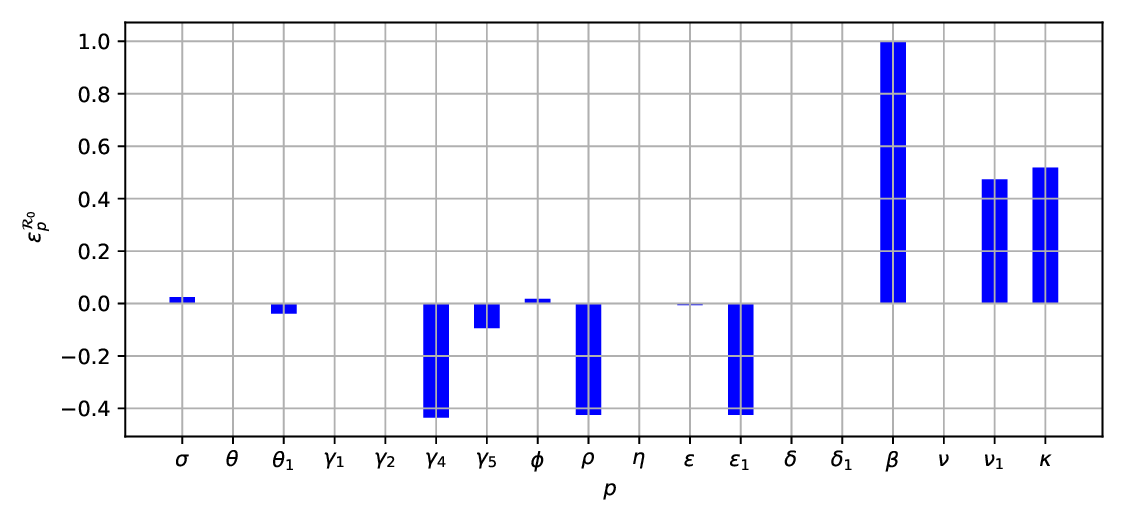}
	\caption{Graphical representation of the sensitivity index of $\mathcal{R}_0$.} \label{fig:sensgraph}
\end{figure}

We now to turn to interpreting the daily infection  data and prevalence data presented in Figures  \ref{fig: prev_vs_daily_infection}. 
An examination of the daily infection data reveals that the peaks of successive waves continued to rise even after the introduction of vaccines, up to the fourth wave (see Figure \ref{fig: sa ddaily infection}).
However, the trend in daily prevalence data shows the opposite; a steady decline over the dame period.
This apparent discrepancy suggests that although the number of new daily infections increased, the overall recovery rate also improved significantly, leading to lower net prevalence. 
Our model supports this interpretation, indicating that recovery rates increased and death rates declined following the rollout of vaccines. 
for complete records of prevalence and mortality data, we refer the reader to the   worldometer database \cite{worldometer}. 

An unexpected, though  not significant,   observation is the positive relation between the vaccination rate parameter $\sigma$ and the basic reproduction number $\mathcal{R}_0.$ One possible explanation for the positive sensitivity index   $ \epsilon^{\mathcal{R}_0}_\sigma$ is that vaccinated individuals, due to incomplete protection against infection, still contribute to transmission, particularly through asymptomatic infections. 
Despite this, vaccination remains critical; our model indicates that the improvement in recovery rates following vaccination led to decline in prevalence and ultimately helped curb the spread of the disease.
Thus, the implementation of the vaccine program was essential for reducing mortality and saving lives.
%
Furthermore, the sensitivity indices   $\epsilon^{\mathcal{R}_0}_{\nu_1}$ and $\epsilon^{\mathcal{R}_0}_\kappa$ highlight another important aspect; vaccinated individuals, feeling relatively safe, may have increased their contact with the public. If infected, such behavior could enhance transmission, contributing to the observed positivity of  $ \epsilon^{\mathcal{R}_0}_\sigma$. This underlines the importance of continued public health messaging, even in highly vaccinated populations, to mitigate behavioral risks.

Based on the findings of this study, we conclude that greater public  and stakeholder awareness of the actual level of protection provided by COVID-19 vaccines at the time could have led to a different trajectory of disease spread.
If individuals  had understood that vaccine effectiveness might be lower than 50\%, rather than assuming it was above 50\%, they may have adopted additional protective measures to safeguard themselves and those around them.
For further concluding remarks, we refer to Section~\ref{conclusion}.

\begin{footnotesize} 
	\begin{table}
		\centering
		\begin{tabular}{>{\stepcounter{rowno}\therowno.}c|cp{2cm}p{2.7cm}p{1.0cm}p{1.3cm}p{1.5cm}p{1.6cm}}
			\hline
			\multicolumn{1}{r}{No.}&Param &  Initial Estimate & Fitting range& fitted value& $\gamma^{\mathcal{R}_0}_{param}$  &$\epsilon^{\mathcal{R}_0}_{param}$  & Source \\
			\hline
			\hline
			&$\Lambda$ &  3538.3 & N/A &N/A & 0& 0 &calculated \\  
			&$\sigma$ &$ 7.9\times10^{-4}$&N/A&N/A&185.9& 0.0242&   calculated \\   
			&$ \mu$ &$2.6433\times10^{-5}$& N/A&N/A & -5603.6&-0.0244&   calculated \\   
			&$\theta$ &0.0267&(0.01,0.03)& 0.03 &-0.0847&-0.0004&  fitted \\   
			&$\theta_1$ &0.0267&$(0.01,\theta)$& 0.01& -23.2&-0.0382& fitted  \\  
			&$\gamma_1 $ &0.0904&(0.0544,0.1167)& 0.1167 &-0.0847& -0.0016& fitted\\  
			&$  \gamma_2 $ &0.0175& (0.0066,0.0313)& 0.0066&-0.3303&-0.0004&fitted \\  
			&$ \gamma_3$ &0.09 & (0.0694,0.0974) & 0.0974 & 0&0 & fitted \\  
			&$\gamma_4$ & 0.095 &(0.0544,0.1167) & 0.1141 &-23.20& -0.4354& fitted\\
			&$  \gamma_5 $ &0.0275 & (0.0066,0.0313) & 0.0233 & -0.0847&-0.0934& fitted \\ %
			&$\omega$ & $1/120 $&N/A&N/A&0&0 &\cite{wajnberg2020robust}\\  
			&$\varphi$ &0.0022&$(0.0011,\omega)$&0.0011&0& 0 & fitted \\  
			&$\phi$ &0.5&$(\eta,1)$& 0.468& 0.2337&0.018&fitted \\  
			&$\rho$ &0.75&(0.3,1)& $0.300$& -8.612&-0.4251& fitted \\  
			&$\eta$ &0.45&(0.3,0.6)&0.30 &-0.0121&-0.0006 & fitted \\  
			&$\epsilon$ &$0.1252$ &(0.1057,0.1472)& 0.1057& -0.3303&-0.0057& fitted  \\   
			&$\epsilon_1$ &$0.1252$&(0.1057,0.1472)& 0.1057& -24.413&-0.4246&fitted \\   
			&$ \delta$ &0.0015&(0.0012,0.0016)&0.0012&-0.3303&-0.00 & fitted \\  
			&$\delta_1$ & 0.0011&$(0,\delta)$& 0.0& -24.414&-0.00  &fitted \\  
			&$\beta$ &0.9&(0,3)&0.1878 &32.36&1.0&fitted \\ 
			&$\nu$ &3.5&(0,6)& 1.0& 0.0124&0.0021&fitted\\ 
			&$\nu_1$ &3.5&(0,6)& 6.0&0.4799& 0.4737&fitted \\ 
			&$\kappa$ &1.25&(1,6)& 6.0& 0.5248&0.5181&fitted \\ 
			\hline
		\end{tabular}
		\caption{Initial estimate of involved parameters, their possible range according to literature, their fitted values and  sensitivity index of $\mathcal{R}_0$ with respect to the involved parameters. }   
		\label{tab: Table3}
	\end{table}
\end{footnotesize}

\section{ Conclusion }\label{conclusion}

In this study, we developed and analyzed a mathematical model   that captures   the imperfect efficacy of COVID-19 vaccines, incorporating key features   such as waning immunity and breakthrough infections to evaluate their impact on pandemic trajectory in South Africa. 
The model integrates multiple dimensions of disease dynamics, including viral transmission, vaccine mediated protection against sever outcomes, public health containment measures, and alignment with WHO and CDC guidelines.
Although the model was calibrated using South African case data, a decision motivated by the authors' regional familiarity and data availability, its structure sufficiently flexible to be adapted to other settings that deployed similar vaccines and implemented comparable distribution strategies.

The model was analyzed both theoretically and numerically.  
The theoretical analysis demonstrated that   when vaccines are imperfect,  the model may exhibit a backward bifurcation even when the basic reproduction number  $\mathcal{R}_0<1$ (see Theorem \ref{Theorem3.33} and Theorem \ref{Theorem3.77}). 
However,  in the absence of a backward bifurcation,  the disease-free equilibrium is shown to be globally asymptotically stable  whenever $\mathcal{R}_0 < 1$ (see Theorem \ref{Theorem3.4}).  
Remark \ref{Theorem3.7}  establishes the existence of at least one locally asymptotically stable endemic equilibrium when   $\mathcal{R}_0>1$ and  close to unity.   
  Furthermore, the positive impact of perfect vaccination in reducing  disease prevalence  is demonstrated in Theorem \ref{Theorem3.8} . 

For numerical analysis, we fitted the South African COVID-19 data to model~\eqref{eq:(2.1)}. 
The parameter estimation results indicate that vaccination significantly improved recovery outcomes, with vaccinated individuals experiencing faster recovery ($\gamma_5>\gamma_2$ ) and lower mortality ($\delta_1<\delta$) compared to their unvaccinated counterparts (see Table \ref{tab: Table3}, rows 7, 10, 18, and 19). These results provide quantitative confirmation that vaccination reduces both disease severity and the risk of death.

The local sensitivity indices of $\mathcal{R}_0$ with respect to the model parameters were computed to identify those with greatest influence. The result indicate that the parameters $\phi,\,\beta,\, \nu_1,\, \kappa,\, \rho,\,\gamma_4, \theta_1$ and $\epsilon_1$ are the most significant. Among these,  $\phi,\,\beta,\, \nu_1$ and $ \kappa$ are positively correlated with  $\mathcal{R}_0$, while   $\rho,\,\gamma_4, \theta_1$ and $\epsilon_1$ are negatively correlated.   
This yields the following proportional relationships:
\[\displaystyle \mathcal{R}_0 \propto \;\beta,\; \nu_1,\;  \;\kappa \;~~~
\mbox{and}\;~~~
\mathcal{R}_0 \propto\frac1\rho,\; \frac{1}{\epsilon_1},\;\frac{1}{\gamma_4}.\]
A graphical representation of these relationships is provided in Figure ~\ref{fig: vaccine effect2}.

%
 Therefore, an effective intervention strategy should aim to reduce  the value of parameters $\beta,\, \nu_1, $ and $\kappa$, while enhancing the  protective  parameters $\rho,\;\gamma_4,\; \theta_1$ and $\epsilon_1.$
The following measures could contribute to achieving these objectives: 
\begin{description}
	\item[Support vaccine research and development:]  investing in scientific research to produce more effective vaccines can help increase the value of $\rho,$ thereby enhancing protection against infection. 
	\item[Promote adherence to public health guidlines:] Encouraging the public to follow WHO recommended non-pharmaceutical interventions, such as avoiding unnecessary social contact even after vaccination can contribute to reducing the transmission rate $\beta.$ 
	\item[Raise awareness of vaccine limitations:] Educating the public about the actual efficacy levels of COVID-19 vaccines can reduce overconfidence among vaccinated individuals. This can lead to a reduction in  $\nu_1$ and $\kappa$, as more cautions behavior among vaccinated individuals lowers the risk of virus transmission. 
	\item[Expand testing and contact tracing:] Increasing the frequency of testing, especially among vaccinated individuals, and implementing effective contact tracing can help identify and isolate infected persons. This would enhance the detection and isolation rates $\theta_1$ and  $\epsilon_1$, ultimately contributing to the reduction of $\mathcal{R}_0$. 
\end{description}

\begin{figure}[h!]
	\centering
	\begin{subfigure}[b]{0.45\textwidth}
		\includegraphics[width=\textwidth]{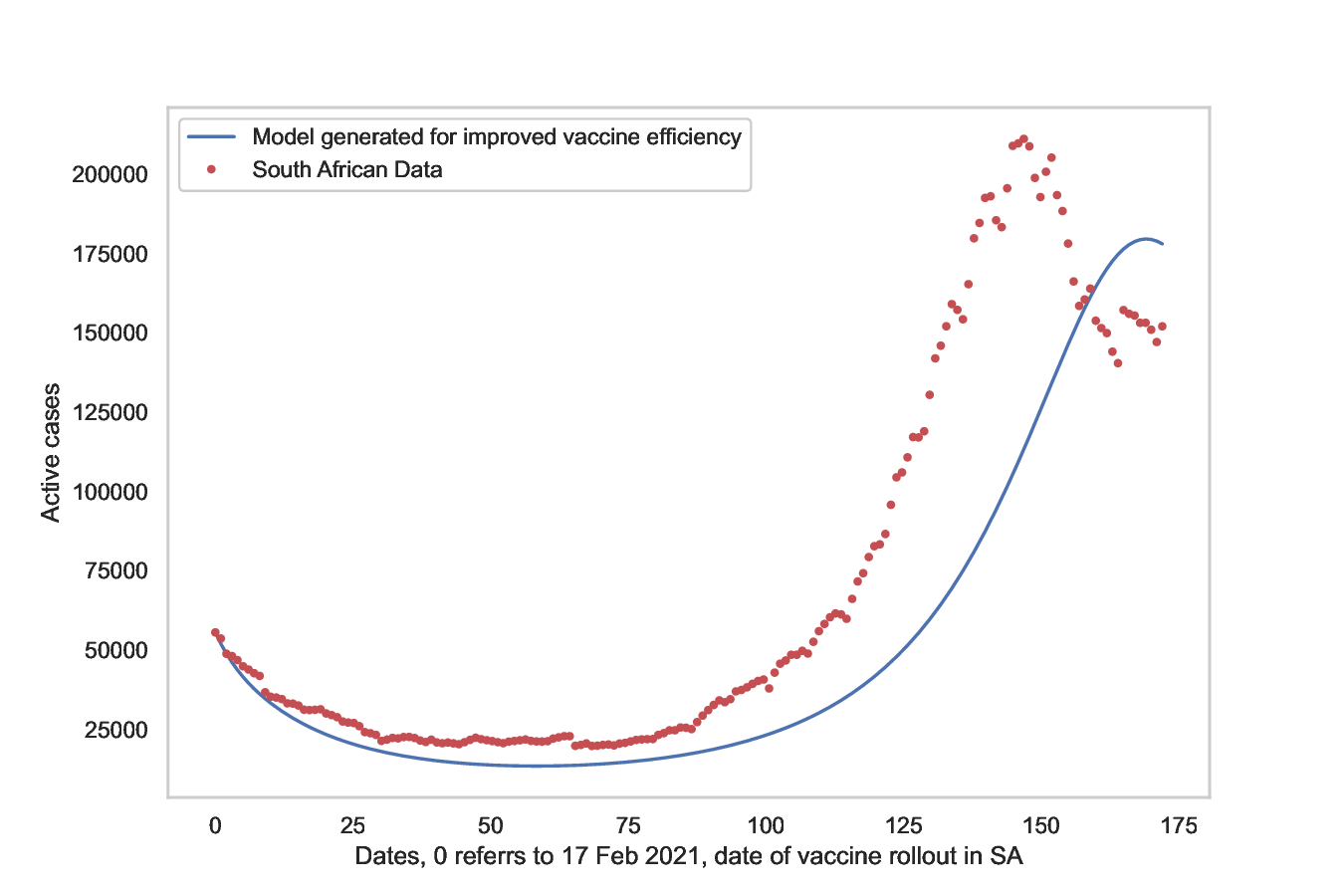}
		\caption{40\% efficacy}
		\label{fig: fitted rho0.25}
	\end{subfigure}
	~
	\begin{subfigure}[b]{0.45\textwidth}
		\includegraphics[width=\textwidth]{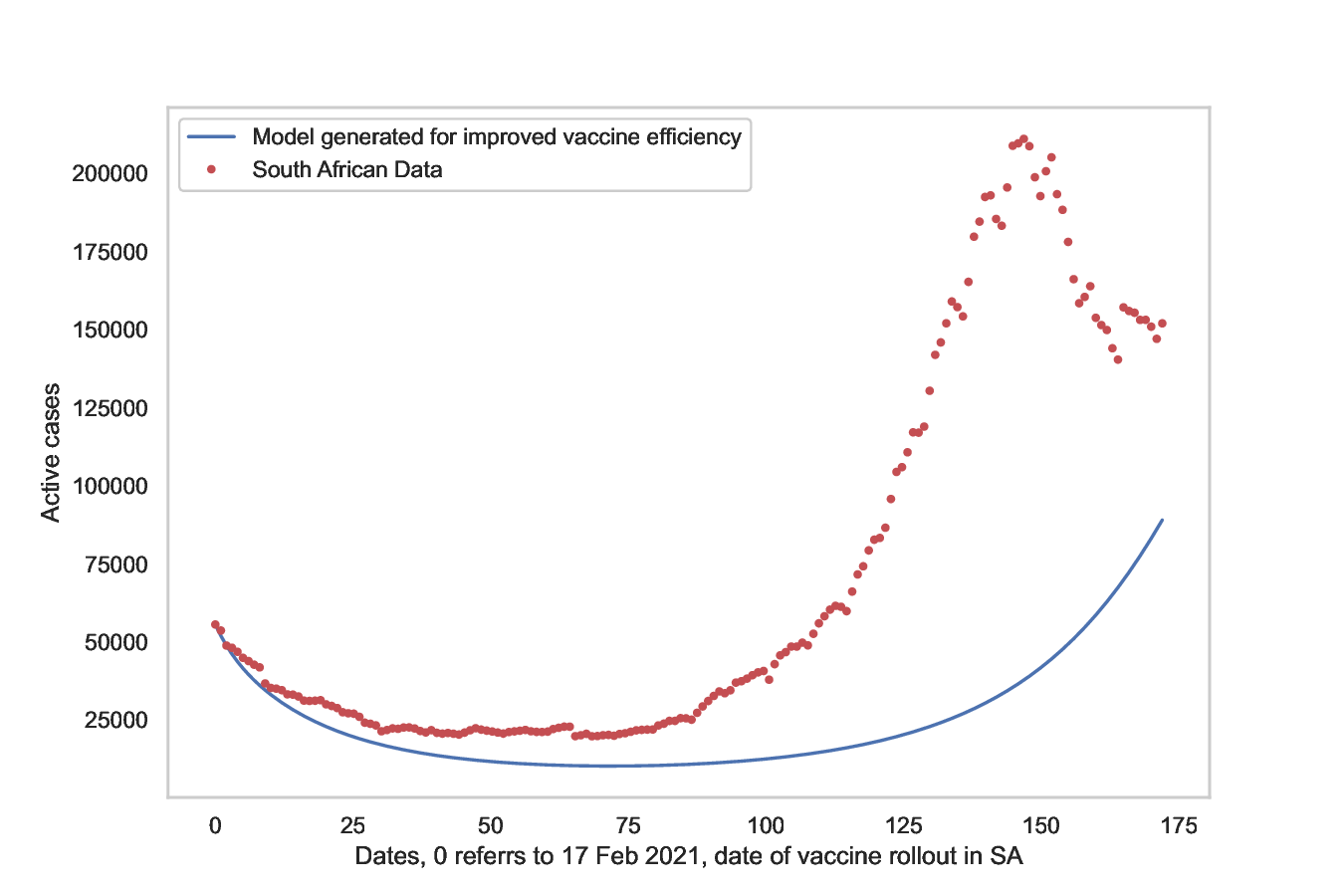}
		\caption{50\% efficacy}
		\label{fig: fitted rho050}
	\end{subfigure}
	
	\begin{subfigure}[b]{0.45\textwidth}
		\includegraphics[width=\textwidth]{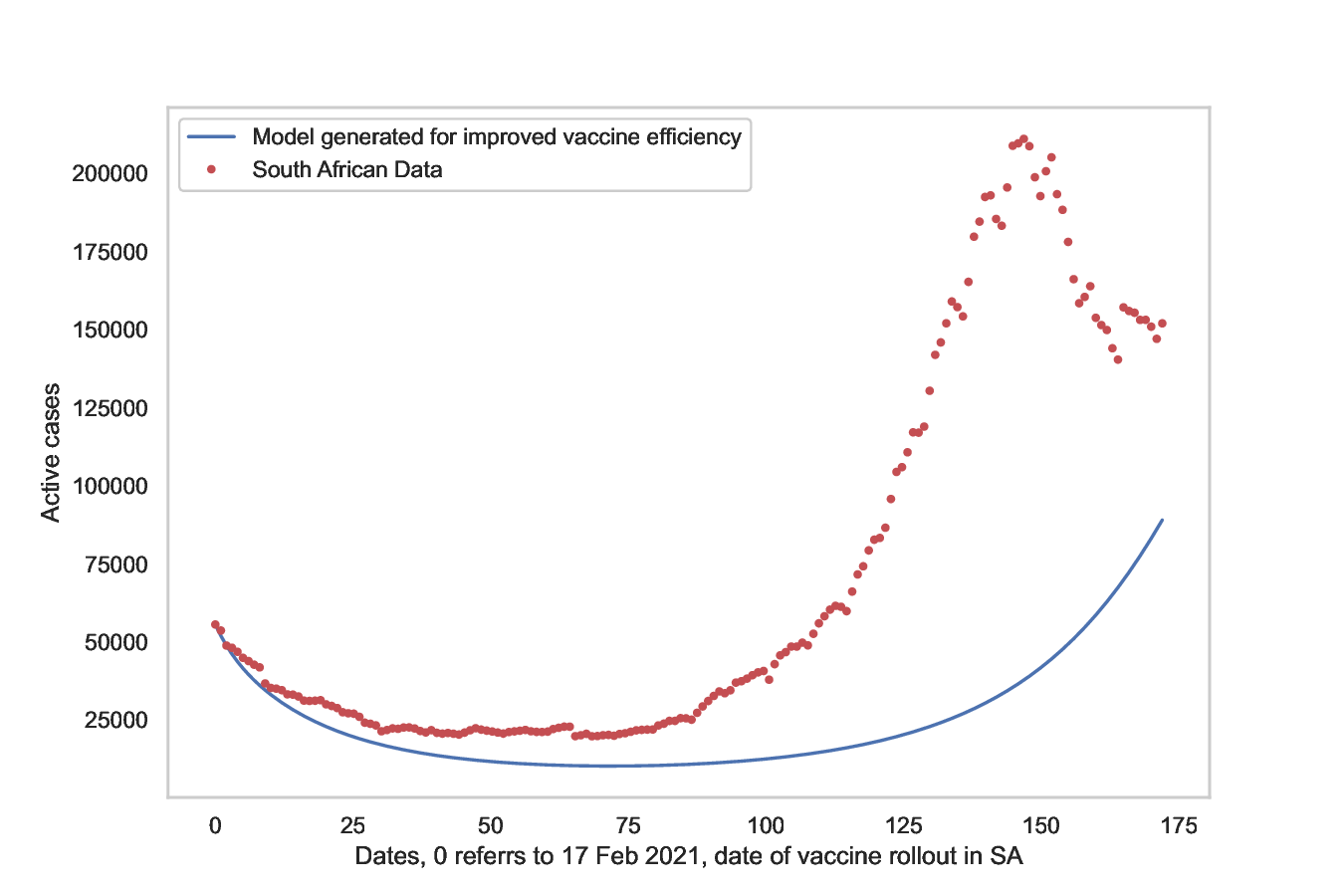}
		\caption{75\% efficacy}
		\label{fig: fitted rho075}
	\end{subfigure}
	\begin{subfigure}[b]{0.45\textwidth}
		\includegraphics[width=\textwidth]{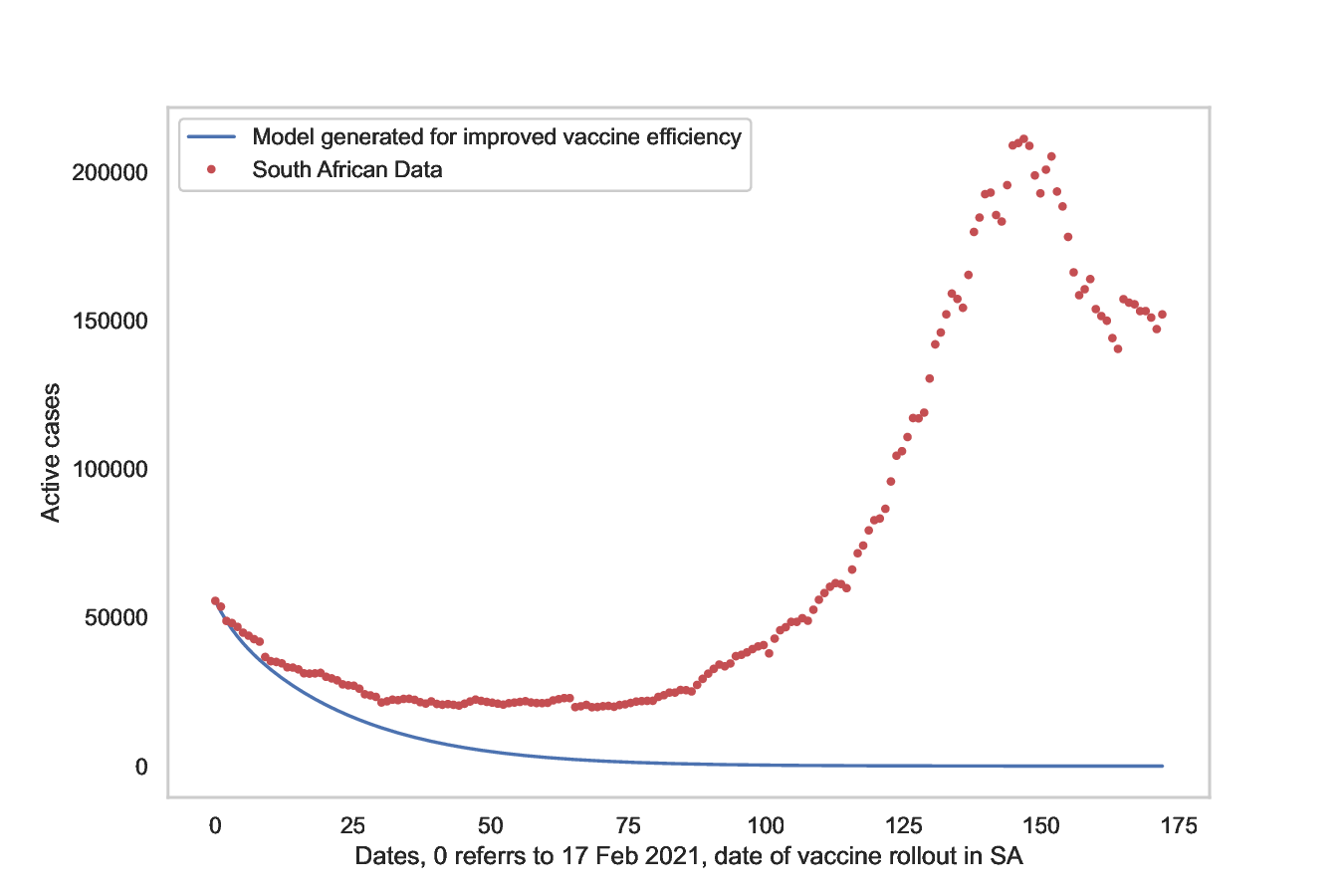}
		\caption{Full protection and no reinfection}
		\label{fig: fitted rhoIdeal}
	\end{subfigure}
	\caption{ Different efficacy levels were considered to demonstrate the role that vaccines could play in flattening the curve and delaying the peak.}\label{fig: vaccine effect}
\end{figure}

\begin{figure}[h!]
	\centering 
	\begin{subfigure}[b]{0.3\textwidth}
		\includegraphics[width=\textwidth]{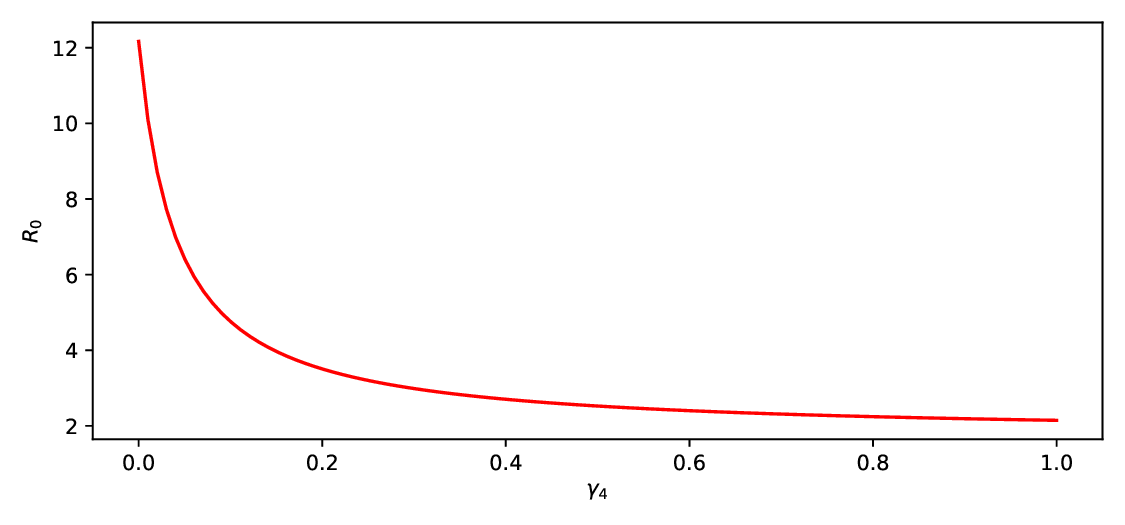}
		\caption{Graphical visualizing the sensitivity of $R_0$ to $\gamma_4$}
		\label{fig: R0gamma4}
	\end{subfigure} \qquad
	\begin{subfigure}[b]{0.3\textwidth}
		\includegraphics[width=\textwidth]{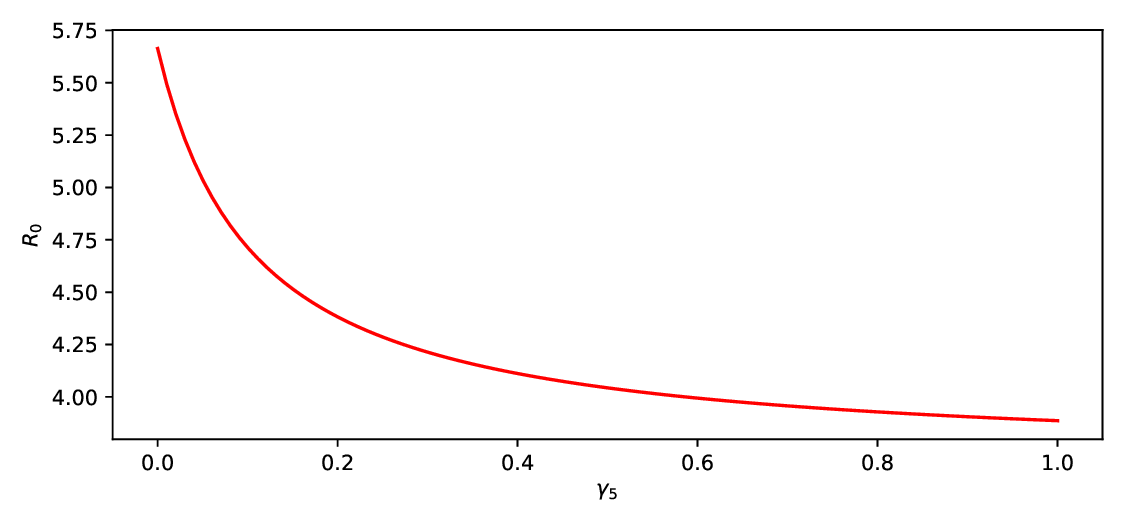}
		\caption{Graphical visualizing the sensitivity of $R_0$ to $\gamma_5$}
		\label{fig: R0gamma5}
	\end{subfigure}\qquad
	~ 
	\begin{subfigure}[b]{0.3\textwidth}
		\includegraphics[width=\textwidth]{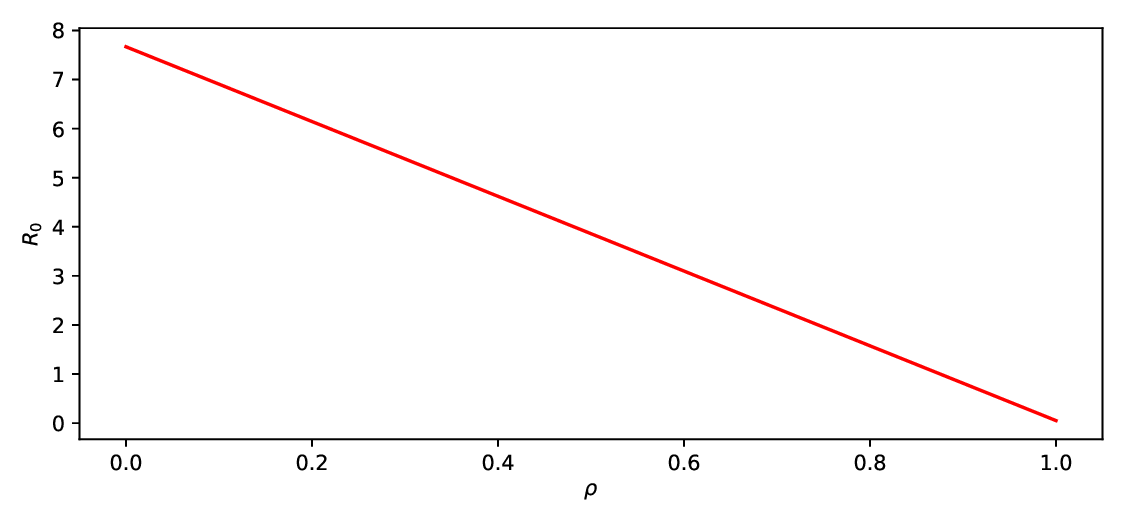}
		\caption{Graphical visualizing the sensitivity of $R_0$ to $\rho$}
		\label{fig: R0rho}
	\end{subfigure}\qquad
	\begin{subfigure}[b]{0.3\textwidth}
		\includegraphics[width=\textwidth]{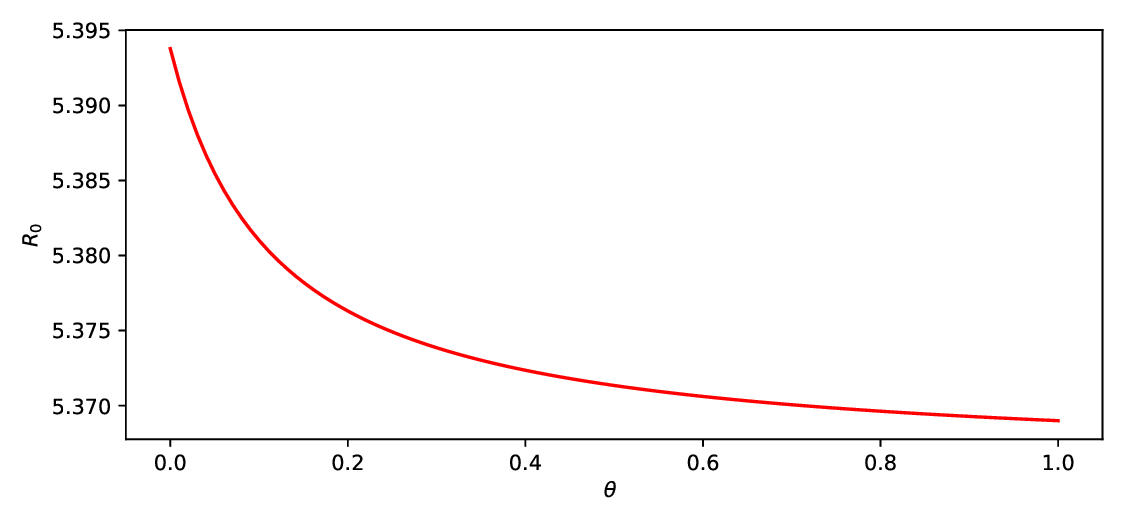}
		\caption{Graphical visualizing the sensitivity of $R_0$ to $\theta$}
		\label{fig: R0theta}
	\end{subfigure}\qquad
	\begin{subfigure}[b]{0.3\textwidth}
		\includegraphics[width=\textwidth]{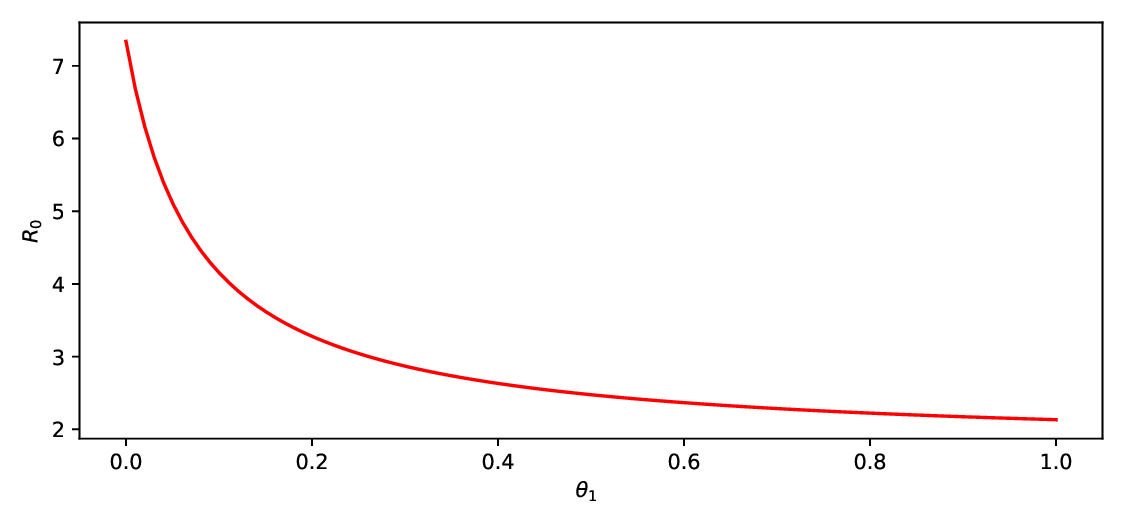}
		\caption{Graphical visualizing the sensitivity of $R_0$ to $\theta_1$}
		\label{fig: R0theta1}
	\end{subfigure}\qquad
	\begin{subfigure}[b]{0.3\textwidth}
		\includegraphics[width=\textwidth]{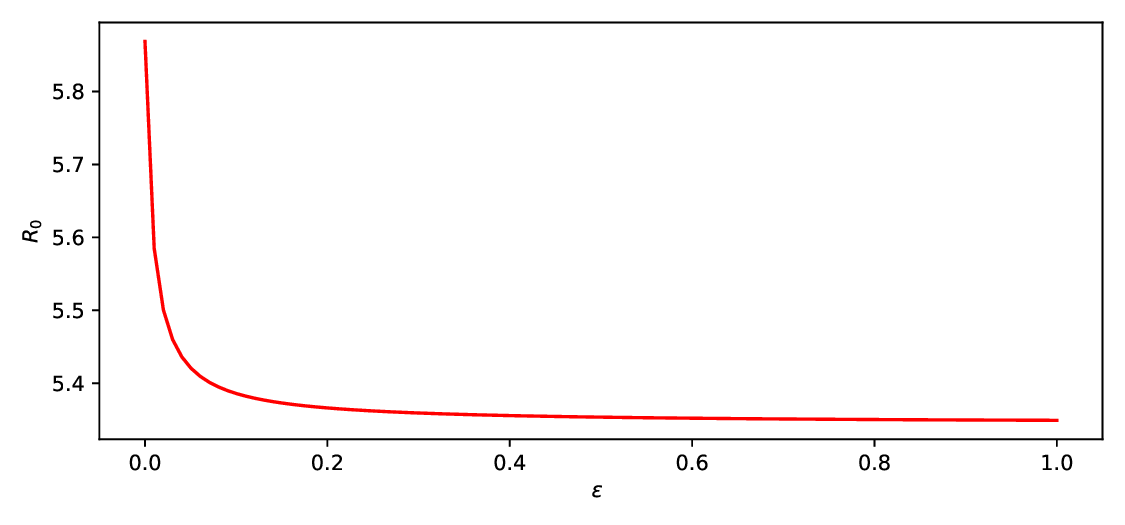}
		\caption{Graphical visualizing the sensitivity of $R_0$ to $\epsilon$}
		\label{fig: R0epsilon}
	\end{subfigure}\qquad
	\begin{subfigure}[b]{0.3\textwidth}
		\includegraphics[width=\textwidth]{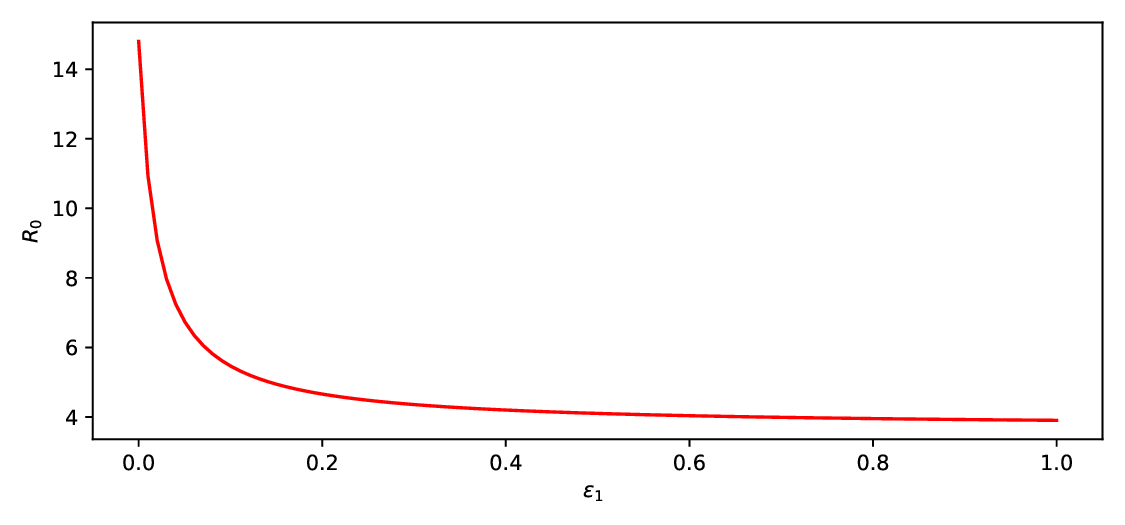}
		\caption{Graphical visualizing the sensitivity of $R_0$ to $\epsilon_1$}
		\label{fig: R0epsilon1}
	\end{subfigure}
	
	\caption{${\cal R}_0$  versus some parameters involved in equation \ref{BRN}.} \label{fig: vaccine effect2}
\end{figure}

In general, based on both  numerical and theoretical   analysis of our model, we draw the following conclusions:
\begin{enumerate} [(1)]
	\item  The result $\delta > \delta_1$ indicates that  vaccines improved recovery rate. This finding aligns with official reports from the WHO and CDC, such as \cite{CDCvaccine,whovaccine}.
	\item The estimate $\rho \approx 0.3$ suggests that that vaccines did not provide  complete protection from COVID-19 infection, highlighting the need for continued research to enhance vaccine efficacy. This observation is also consistent with findings from the   NCIRD report \cite{NCIRD}.
	\item The result $\gamma_5>\gamma_2$ shows that vaccinated individuals recover more quickly and have a reduced risk of mortality. This outcome is also supported by WHO and CDC data \cite{CDCvaccine,whovaccine}.
\end{enumerate}

Based on the results of this work, we   recommend  that  public health policy makers and other stakeholders   consider the following points in their strategic planning: 
\begin{enumerate}[(1)]
	\item Encourage the public to reduce mobility and social interactions when feasible to limit contact with potential infectious individuals, especially in the event of future outbreaks. This aligns with the fact that    $\mathcal{R}_0\propto\beta,\;\mathcal{R}_0\propto\nu, \;\mathcal{R}_0\propto\nu_1\; \text{and }\; \mathcal{R}_0\propto\kappa$.  
	\item   Provide increased support for research initiatives aimed at improving  the effectiveness of   COVID-19   vaccines, particularly with respect to protection against infection. 
	\item   Ensure the availability of accessible and  affordable   testing    centers. Increased testing enhances the ability to identify and isolate infected individuals, thereby increasing the parameters $\theta$ and $\theta_1$.  
	\item    Establish isolation facilities for individuals who are unable to self-isolate, such as those who are homeless or living in overcrowded households. Under ideal conditions, isolated individuals do not contribute to the spread of infection  (see equation \eqref{BRN}).
\end{enumerate}

Finally, this study offers several avenues for future research. These include the identification of cost-effective intervention   strategies that optimize disease control  while minimizing economic burden, as well as the development of stochastic models to incorporate uncertainty in transmission dynamics and intervention effectiveness.   %

\section*{Acknowledgment}
 The authors gratefully acknowledge the anonymous reviewers and the Associate Editor for their insightful comments and suggestions, which have significantly enhanced this work.
Part of this work was presented at the workshop ``On Research Trend in Mathematical Modeling and Analysis in Life Sciences'' held at Tshikwalo Game Lodge, Dinokeng Game Reserve, Pretoria, South Africa and  at the ``Second annual One Health Symposium'' held at 
University of Pretoria, Pretoria, South Africa. The   first author acknowledges the following grants for their financial support;
\begin{itemize}
	\item  DSI-NRF Center of Excellence in Mathematical and Statistical Sciences (CoE-Mass) ref. No. 2022-003-21F-trends.
	\item  SARCHI chair in Mathematical Models and Methods in Bio-engineering and Bio-sciences (${\rm M^3B^2}$).
\end{itemize}


\appendix  

\section{Proof of Theorem \ref{Theorem3.1}} \label{Appendix A}

\begin{proof}

	We want to prove that for non-negative initial condition, at all time $t\geq0$, the system
	(\ref{eq:(2.1)}) has a unique non-negative solution which is contained in $\Omega$. The proof is done in three steps: Firstly we show the non-negativity of the solutions for any non-negative initial data; secondly we establish the  boundedness of the solution and finally we establish uniqueness of the solution.

	To prove the non-negativity, we use the method of contradiction as it is in  \cite{Busenberg_1993,Terefe_2020}.  
	Without loss of generality, we may assume that the trajectory of $R$ will pass to the negative region before others,   
	i.e., we consider the trajectory of $R$ crosses to the region    $R<0$ at some positive time $t_1$, such that 
	\begin{align}
		&R(t_1)=0,\;R'(t_1)<0,\; \mbox{and}\;\nonumber \\  &A(t_1)>0,\;I(t_1)>0,A_1(t_1)>0,I_1(t_1)>0,\;~Q(t_1)>0. \label{eq: Rpassumption}
	\end{align}
	
	Then, from   equation \eqref{eq:(2.1)R} we have,  %
	\begin{align}
		R'(t_1)=\gamma_1 A(t_1)+\gamma_2I(t_1)+\gamma_3Q(t_1)+\gamma_4A_1(t_1)+\gamma_5I_1(t_1). \label{eq: Rpositive}
	\end{align}
	Observe that, due to   assumption   \eqref{eq: Rpassumption}, the left hand side of \eqref{eq: Rpositive} is negative while the right hand side  is positive, which is a contradiction.  
	Hence, $R(t)$ remains non-negative for all $t\geq0$. 
	From   equation \eqref{eq:(2.1)S} we have
	\begin{eqnarray}
		S'(t)= \Lambda- [\lambda+\sigma+\mu]S +\varphi R \geq - [\lambda+(\sigma+\mu)]S. \label{eq Spositive}
	\end{eqnarray}
	Applying simple calculus techniques to \eqref{eq Spositive}, we obtain
	\begin{eqnarray*}
		S(t)\geq S(0)\exp\left(-\int_0^t\left(\lambda(u)+\sigma+\mu\right)du\right)\geq0.
	\end{eqnarray*}
	Thus,  $S(t)$ remains non-negative for all $t\geq0$. Similarly, from   equation \eqref{eq:(2.1)V} we have
	\begin{eqnarray}
		V'(t)= \sigma S-[(1-\rho)\lambda  +\mu ]V+(\omega-\varphi)R\geq-[(1-\rho)\lambda  +\mu ]V ,  \label{eq:Vpositive}
	\end{eqnarray}
	which yields
	\begin{eqnarray*}
		V(t)\geq V(0)\exp\left(-\int_0^t\left((1-\rho)\lambda(u)  +\mu\right)du \right)\geq0.
	\end{eqnarray*}
	Hence, $V(t)$ also remains non-negative for any $t\geq0$. 
	
	To show the non-negativity of the variables $A,\, A_1,\, I,\, I_1 $ and $Q$ one can follow a procedure similar to the one used to show the non-negativity of $R.$

	To proof the boundedness of the system, we use principle of conservation. From \eqref{eq: total population} and \eqref{eq:(2.1)},  
	we obtain 
	\begin{eqnarray} \label{Gron}
		N'(t)=\Lambda-\mu N-\delta(I+Q)-\delta_1I_1\leq \Lambda-\mu N.
	\end{eqnarray}
	For an initial population $N_0$, implementing Gronwall's inequality on \eqref{Gron} gives 
	\begin{align}
		N(t)\leq\Lambda/\mu+(N_0-\Lambda/\mu)\exp(-\mu t) <\infty .\label{eq: Nbonded}  
	\end{align}
	Hence, the solution of the model is bounded for every time $t\geq0$.

	Finally, the uniqueness follows from Steps 1 and 2 and %
	Theorem 2.1.5 of \cite{Stuart_1998}. Thus we are guaranteed that any solution of \eqref{eq:(2.1)} is non-negative and bounded for $t\geq0$. Thus, the model equation \eqref{eq:(2.1)} is a dynamical system on $\Omega$. This completes the proof of Theorem \ref{Theorem3.1}.
	
\end{proof}
\begin{remark}
	Equation\eqref{eq: Nbonded} tells us that   the total population at any given time $N(t)$ remains bounded. In fact 
	\begin{align*}
		\lim_{t\to\infty} \left( \Lambda/\mu+(N_0-\Lambda/\mu)\exp(-\mu t) \right)= \frac{\Lambda}{\mu}.
	\end{align*}
	Thus, the set   
	{\small
		\begin{align*}
			\tilde{\Omega} = \Big\{(S,V, A,I,A_1,I_1, Q,  R) \in \mathbb{R}_{+}^{8}:
			0\leq S+V+A+I+A_1+I_1+Q+R=N  \leq \frac{\Lambda}{\mu} \Big\} \subset\Omega
		\end{align*}
	}
	is an attractor set of the system \eqref{eq:(2.1)}.   
\end{remark}
\section{Calculation of the Basic Reproduction Number}\label{BRN calculation}
The basic reproduction number refers the average number of secondary cases produced in a completely susceptible population by an  infectious individual during his/her entire infectious period \cite{van_2008}.    
We compute $\mathcal{R}_0$ by using the method of Next Generation Matrix,  see  \cite{Castilo_2002, Diekmann_2000, van_2008, Shuai_2013}.

For the model under consideration, we denote infected classes by $\mathcal{A}$  and define vector valued functions $\mathcal{F}:\mathcal{A}\to \mathbb{R}^5 $ and $\mathcal{U}: \mathcal{A}\to \mathbb{R}^5$ by

\begin{align*}
	\mathcal{F}(X)= \begin{pmatrix}
		\eta\lambda S\\
		(1-\eta)\lambda S\\
		\phi(1-\rho)\lambda V\\
		(1-\phi)(1-\rho)\lambda V\\
		0
	\end{pmatrix}  \quad \text{ and }\quad \mathcal{U}(X)=\begin{pmatrix}
		k_1A\\
		k_2 I\\
		k_3 A_1\\
		k_4 I_1\\
		-\theta A-\theta_1 A_1-\epsilon I-\epsilon_1 I_1+k_5Q
	\end{pmatrix}
\end{align*}
where \begin{align*}
	\mathcal{A} = \big\{(A, I,A_1, I_1,Q): (S,V,A, I, A_1, I_1,Q, R)\in \Omega  \big\}.
\end{align*}
The function  $\mathcal{F}$ represents the rate of appearance of new infection and $\mathcal{U}$ denotes the rate of transfer of individuals among the infective classes, respectively, where
$$k_1=\theta +\gamma_1+\mu,\;k_2=\epsilon+\gamma_2+\delta+\mu,\;k_3=\theta_1+\gamma_4+\mu,\;k_4=\epsilon_1+\gamma_5+\delta_1+\mu$$ and $$k_5=\gamma_3+\delta+\mu.$$

The next generation matrix is given by
\begin{equation} \label{$(3.6)$}
	\mathcal{K}=J_\mathcal{F}J_\mathcal{U}^{-1},
\end{equation}
where
\begin{align}
	J_\mathcal{F}= \begin{pmatrix}
		B_1\nu&B_1&B_1\nu_1 &B_1\kappa&0\\
		B_2\nu&B_2&B_2\nu_1&B_2\kappa&0 \\
		B_3\nu&B_3&B_3\nu_1&B_3\kappa&0 \\
		B_4\nu&B_4&B_4\nu_1&B_4\kappa&0 \\
		0&0&0&0&0
	\end{pmatrix}  \quad \mbox{and} \quad
	J_\mathcal{U} = \begin{pmatrix}
		k_1&0&0&0&0\\
		0&k_2&0&0&0\\
		0&0&k_3&0&0\\
		0&0&0&k_4&0\\
		-\theta&-\theta_1 & -\epsilon&-\epsilon_1&k_5
	\end{pmatrix} \label{mat Jf}
\end{align}
are the Jacobian matrices  of $\mathcal{F}$ and $\mathcal{U}$ at $E_0$, respectively with \[B_1=\frac{\eta \beta \mu}{\sigma+\mu},\;B_2=\frac{(1-\eta) \beta \mu}{\sigma+\mu},\;B_3=\frac{\phi(1-\rho) \beta \sigma}{\sigma+\mu}\; \mbox{and} \;B_4=\frac{(1-\phi)(1-\rho) \beta \sigma}{\sigma+\mu}.\] Notice that 

\begin{align}
	J_\mathcal{U}^{-1}=  \begin{pmatrix}
		\frac{1}{k_1}&0&0&0&0\\ 0&\frac{1}{k_2}&0&0&0\\0&0&\frac{1}{k_3}&0&0\\0&0&0&\frac{1}{k_4}&0\\ \frac{\theta}{k_1k_5}&\frac{\theta_1}{k_2k_5}&\frac{\epsilon}{k_3k_5}&\frac{\epsilon_1}{k_4k_5}&\frac{1}{k_5}
	\end{pmatrix} .\label{mat Ju-1}
\end{align}

We now  combine equations \eqref{mat Jf} and \eqref{mat Ju-1}, to get
\begin{align}
	\mathcal{K}= \begin{pmatrix}
		\frac{B_1\nu}{k_1}&\frac{B_1}{k_2}&\frac{B_1\nu_1}{k_3} &\frac{B_1\kappa}{k_4}&0\\
		\frac{B_2\nu}{k_1}&\frac{B_2}{k_2}&\frac{B_2\nu_1}{k_3}&\frac{B_2\kappa}{k_4}&0 \\
		\frac{B_3\nu}{k_1}&\frac{B_3}{k_2}&\frac{B_3\nu_1}{k_3}&\frac{B_3\kappa}{k_4}&0 \\
		\frac{B_4\nu}{k_1}&\frac{B_4}{k_2}&\frac{B_4\nu_1}{k_3}&\frac{B_4\kappa}{k_4}&0 \\
		0&0&0&0&0
	\end{pmatrix}.
\end{align}
Thus, we have  
\begin{eqnarray} 
	\mathcal{R}_0= \mathcal{R}_A + \mathcal{R}_I + \mathcal{R}_{A_1} + \mathcal{R}_{I_1}  \label{BRNA}
\end{eqnarray}
where $$\displaystyle \mathcal{R}_A = \frac{\nu B_1 }{k_1},\,~  \mathcal{R}_I = \frac{B_2}{k_2},\,~ \mathcal{R}_{A_1} = \frac{\nu_1 B_3}{k_3},\,~ \mathcal{R}_{I_1} = \frac{\kappa B_4}{k_4}.$$

				\section{Proof of Theorem \ref{Theorem3.77}} \label{Appendix BB}
					The proof is based on the Center manifold theory and Theorem 4.1 from \cite{Castilo_2004}. 
				For this we introduce a change of variables by setting , $x_1=S,\;x_2=V,\;x_3=A,\;x_4=I,\;x_5=A_1,\;x_6=I_1,\;x_7=Q$ and $x_8=R$ %
				and use the new variables to write model \eqref{eq:(2.1)} as%
				\[X'(t)=f(X)=(f_1,f_2,f_3,f_4,f_5,f_6,f_7,f_8)^T,\]  
				with $X=(x_1,x_2,x_3,x_4,x_5,x_6,x_7,x_8)^T$, and %
					\begin{equation} \label{eq:2.1a1} 
						\begin{aligned}
							f_1:=x_1'(t) &= \Lambda- [\lambda+k_0]x_1 +\varphi x_8,\\
							f_2:=x_2'(t)  &= \sigma x_1-[(1-\rho)\lambda  +\mu ]x_2+(\omega-\varphi)x_8,\\
							f_3:=x_3'(t)  &= \eta\lambda x_1 -k_1 x_3,\\
							f_4:=x_4'(t)  &= (1-\eta)\lambda x_1 -k_2x_4,\\
							f_5:=x_5'(t) &= \phi(1-\rho)\lambda x_2 -k_3 x_5,\\
							f_6:=x_6'(t)  &= (1-\phi)(1-\rho)\lambda x_2-k_4x_6,\\
							f_7:=x_7'(t)  &= \theta x_3+\theta_1x_5 +\epsilon x_4+\epsilon_1x_5 - k_5x_7,\\
							f_8:=x_8'(t)  &= \gamma_1 x_3 + \gamma_2x_4+\gamma_3 x_7+\gamma_4x_5+\gamma_5x_6 - k_6 x_8,\\
						\end{aligned}
					\end{equation}
					where   $k_0=\sigma + \mu,\;
					k_1 = \theta + \gamma_1 + \mu,\; 
					k_2 = \epsilon + \gamma_2 + \delta + \mu,\;
					k_3 = \theta_1 + \gamma_4 + \mu,\;
					k_4 = \epsilon_1 + \gamma_5 + \delta_1 + \mu,\;
					k_5=\gamma_3 + \delta + \mu,\; $ and 
					$k_6=\omega + \mu$.
					Then, the force of infection $\lambda$ in \eqref{FF} and the disease-free equilibrium $E_0$ in \eqref{DFE}  in terms of the new variables are given by 
					\begin{eqnarray}\label{FF1}
						\lambda =\beta \frac{x_4 +\nu x_3+\nu_1x_5 +\kappa x_6 }{N },
					\end{eqnarray} 
					and 
					\begin{eqnarray}\label{DFE1}
						E_0=\left(x_1,x_2,x_3,x_4,x_5,x_6,x_7,x_8\right)=\left(\frac{\Lambda}{\sigma+\mu}, \frac{\sigma\Lambda}{\mu(\sigma+\mu)},0,0,0,0,0,0\right),
					\end{eqnarray}
					respectively.

					Next, we calculate the  Jacobean matrix of  the system (\ref{eq:2.1a1}) at  \eqref{DFE1} and we get
					{\scriptsize
						\begin{equation}
							J(E_0) = \begin{pmatrix} \label{JM}
								-k_0 & 0 & - \nu m_1 & -m_1 & - \nu_1 m_1 & -\kappa m_1 & 0 & \varphi \\
								\sigma & -\mu & - (1-\rho)\nu m_2 & - (1-\rho)V_0 m_2 & - (1-\rho)\nu_1 m_2 & - (1-\rho)\kappa m_2 & 0 & \omega - \varphi \\
								0 & 0 & \eta  \nu m_1 - k_1 & \eta m_1 & \eta \nu_1 m_1 & \eta  \kappa m_1 & 0 & 0 \\
								0 & 0 & (1-\eta) \nu m_1& (1-\eta)m_1- k_2 & (1-\eta) \nu_1 m_1 & (1-\eta) \kappa m_1 & 0 & 0 \\
								0 & 0 & \phi(1-\rho)\nu m_2 & \phi(1-\rho)m_2 & \phi(1-\rho)\nu_1 m_2- k_3 & \phi(1-\rho) \kappa m_2 & 0 & 0 \\
								0 & 0 & (1-\phi)(1-\rho)\nu m_2& (1-\phi)(1-\rho) m_2 & (1-\phi)(1-\rho)\nu_1 m_2 & (1-\phi)(1-\rho)\kappa m_2 - k_4 & 0 & 0 \\
								0 & 0 & \theta & \epsilon & \theta_1 & \epsilon_1 & -k_5 & 0 \\
								0 & 0 & \gamma_1 & \gamma_2 & \gamma_4 & \gamma_5 & \gamma_3 & -k_6
							\end{pmatrix},
					\end{equation}}
					where \(m_1= \frac{\mu \beta}{\sigma+\mu}\) and \(m_2=\frac{\sigma \beta}{\sigma+\mu}\).
					We now consider  $\beta^*:=\beta$ as a bifurcation parameter at $\mathcal{R}_0 = 1$. Solving for $\beta$  from \eqref{BRNA} at $\mathcal{R}_0=1$ gives
					\begin{eqnarray}\label{par}
						\beta = \frac{k_0k_1k_2k_3k_4}{\nu \eta\mu k_2k_3k_4+(1-\eta)\mu k_1k_3k_4+\nu\phi(1-\rho)\sigma k_1k_2k_4+\kappa(1-\phi)(1-\rho)\sigma k_1k_2k_3}.
					\end{eqnarray}
						At $\mathcal{R}_0=1$  the matrix $J(E_0)$ has a simple zero eigenvalue and all other eigenvalues have negative real parts. This can be verified by substituting \eqref{par} into \eqref{JM} and follow the usual way of calculating the eigenvalue of a matrix.%
					By using the notation in Castilo-Chavez and Song (2004), the following computations are carried out.
					Borrowing the notation of \cite{Castilo_2004}, 
					let $w=(w_1,w_2,w_3,w_4,w_5,w_6,w_7,w_8)$ and  $v=(v_1,v_2,v_3,v_4,v_5,v_6,v_7,v_8)^T$ respectively be the left and right eigenvector   associated with the zero eigenvalue of $J(E_0)$ such that $v\cdot w=1$. Therefore, $w$ and $v$ satisfy the system 
					\begin{eqnarray}\label{C5}
						wJ(E_0)=0\;~\mbox{ and}\;~ J(E_0)v=0.
					\end{eqnarray}
					Solving for $v$ and $w$ from \eqref{C5} yields
					%
					\[v_1=\alpha_1 v_3,\; v_2= \alpha_2 v_3,\;~v_3=v_3>0,\;~
					v_4 = \alpha_4 v_3,\;~
					v_5 = \alpha_5 v_3,\;~
					v_6 = \alpha_6 v_3,\;~v_7=\alpha_7 v_3,\;~v_8=\alpha_8 v_3,
					\]
					where 
					
					\begin{align*}
						\alpha_1 &= \frac{-\nu m_1 - m_1 \alpha_4 - \nu_1 m_1 \alpha_5 - \kappa m_1 \alpha_6 + \varphi \alpha_8}{k_0}, \\
						\alpha_2& = \frac{\sigma \alpha_1 - (1-\rho)\nu m_2 - (1-\rho)V_0 m_2 \alpha_4 - (1-\rho)\nu_1 m_2 \alpha_5 - (1-\rho)\kappa m_2 \alpha_6 + (\omega - \varphi) \alpha_8}{\mu},\\
						\alpha_4& = \frac{-(\eta \nu m_1 - k_1) - \eta \nu_1 m_1 \alpha_5 - \eta \kappa m_1 \alpha_6}{\eta m_1}, \\
						\alpha_5& = \frac{(1-\rho)\phi m_2 \left(\frac{\eta \nu m_1 - k_1}{\eta m_1} - \nu\right)}{\phi(1-\rho)\nu_1 m_2 \left(1 - \frac{m_2}{m_1}\right) - k_3}, \\
						\alpha_6& = \frac{(1-\rho)(1-\phi) m_2 \left(\nu - \frac{\eta \nu m_1 - k_1}{\eta m_1}\right)}{(1-\phi)(1-\rho)\kappa m_2 - k_4}.\\
						\alpha_7& = \frac{\theta + \epsilon \alpha_4 + \theta_1 \alpha_5 + \epsilon_1 \alpha_6}{k_5}, \\
						\alpha_8 &= \frac{\gamma_1 + \gamma_2 \alpha_4 + \gamma_4 \alpha_5 + \gamma_5 \alpha_6 + \gamma_3 \alpha_7}{k_6}.
					\end{align*}
					and
					\begin{align*}
						w_1&=w_2=0,\quad w_3=w_3>0,\\
						w_4& = \left[\frac{k_1 - \eta \nu m_1}{(1-\eta)\nu m_1}\right]w_3 - \left[\frac{\phi(1-\rho)\nu m_2}{(1-\eta)\nu m_1}\right]w_5 - \left[\frac{(1-\phi)(1-\rho)\nu m_2}{(1-\eta)\nu m_1}\right]w_6\\
						w_5 &= \left[\frac{\eta \nu_1 m_1 \left(k_4 - (1-\phi)(1-\rho)\kappa m_2\right) - (1-\eta)\nu_1 m_1 \mathcal{C}}{\mathcal{D}}\right]w_3\\
						w_6 &= \left[\frac{\eta \kappa m_1 \left(k_3 - \phi(1-\rho)\nu_1 m_2\right) + (1-\eta)\kappa m_1 \mathcal{C}}{\mathcal{D}}\right]w_3 \\
						w_7&=w_8=0\\
						\mathcal{C} &= \frac{(\eta \nu m_1 - k_1)(1-\rho)\kappa \nu_1 m_2}{(1-\eta)\nu m_1}  \\
						\mathcal{D} &= k_3 k_4 - (1-\rho)m_2\left(1 - \frac{m_2}{\nu m_1}\right)\left[\phi \nu_1 k_4 + (1-\phi)\kappa k_3\right].
					\end{align*}
					Therefore, the bifurcation coefficients $a$ and $b$ as defined in \cite{Castilo_2004} at the disease-free equilibrium \eqref{DFE1}  are given by 
						\begin{eqnarray}\label{C6}
							a& =& \sum_{k,i,j=1}^8 w_k v_i v_j \frac{\partial^2 f_k}{\partial x_i \partial x_j} ,\\
							&=&w_3 \sum_{i,j=3}^6  v_i v_j \frac{\partial^2 f_3}{\partial x_i \partial x_j} +w_4 \sum_{i,j=3}^6  v_i v_j\frac{\partial^2 f_4}{\partial x_i \partial x_j} \nonumber\\
							&& +w_5 \sum_{i,j=3}^6  v_i v_j\frac{\partial^2 f_5}{\partial x_i \partial x_j}+w_6 \sum_{i,j=3}^6  v_i v_j\frac{\partial^2 f_6}{\partial x_i \partial x_j},\nonumber \\
							&=& 2m_1 \Big[ \eta \left( w_3 v_3 v_4 \nu + w_3 v_3 v_5 \nu \nu_1 + w_3 v_3 v_6 \nu \kappa \right) \nonumber \\
							&&+ (1-\eta)  \left( w_4 v_3 v_4 \nu + w_4 v_3 v_5 \nu \nu_1 + w_4 v_3 v_6 \nu \kappa \right) \Big]\nonumber \\
							&& +2m_2\Big[ \phi(1-\rho) \left( w_5 v_3 v_4 \nu + w_5 v_3 v_5 \nu \nu_1 + w_5 v_3 v_6 \nu \kappa \right) \nonumber \\
							&& + (1-\phi)(1-\rho)  \left( w_6 v_3 v_4 \nu + w_6 v_3 v_5 \nu \nu_1 + w_6 v_3 v_6 \nu \kappa \right) \Big].\nonumber\\
							&=&2\nu\Big[m_1\Big(\eta w_3 +(1-\eta)w_4\Big)+m_2\Big(\phi w_5 + (1-\phi)w_6\Big)\Big]v_3\Big(v_4+\nu_1 v_5+\kappa v_6\Big), \nonumber
						\end{eqnarray}
					and
					\begin{eqnarray}
						b& =& \sum_{k,i=1}^8 w_k v_i \frac{\partial^2 f_k}{\partial x_i \partial \beta} (E_0,\beta^*), \\
						& =&w_3 \sum_{i=1}^8  v_i \frac{\partial^2 f_3}{\partial x_i \partial \beta} +w_4 \sum_{i=1}^8  v_i \frac{\partial^2 f_4}{\partial x_i \partial \beta} +w_5 \sum_{i=1}^8  v_i \frac{\partial^2 f_5}{\partial x_i \partial \beta}+w_6 \sum_{i=1}^8  v_i \frac{\partial^2 f_6}{\partial x_i \partial \beta},\nonumber\\
						&=&\frac{\eta \mu}{\sigma+\mu}\Big(v_3w_3\nu+v_4w_3+v_5w_3\nu_1+v_6w_3\kappa\Big)+\frac{(1-\eta) \mu}{\sigma+\mu}\Big(v_3w_4\nu+v_4w_4+v_5w_4\nu_1+v_6w_4\kappa\Big)\nonumber\\
						&+&\frac{\phi(1-\rho)\sigma}{\sigma+\mu}\Big(v_3w_5\nu+v_4w_5+v_5w_5\nu_1+v_6w_5\kappa\Big)+\frac{(1-\phi)(1-\rho) \sigma}{\sigma+\mu}\Big(v_3w_6\nu+v_4w_6+v_5w_6\nu_1+v_6w_6\kappa\Big), \nonumber\\
						&=&\Big[\frac{\mu}{\sigma+\mu}\Big(\eta w_3+(1-\eta)w_4\Big)+\frac{(1-\rho)\sigma}{\sigma+\mu}\Big(\phi w_5+(1-\phi)w_6\Big)\Big]\Big(\nu v_3+v_4+\nu_1 v_5+\kappa v_6\Big). \nonumber\\
						&=& \frac{\mu(\eta + (1-\eta)\alpha_4) + (1-\rho)\sigma(\phi \alpha_5 + (1-\phi)\alpha_6)}{\sigma + \mu} \cdot \Big(\nu + \alpha_4 + \nu_1 \alpha_5 + \kappa \alpha_6\Big)v_3w_3 
					\end{eqnarray}
					where $v_3>0$ and $w_3>0$.
					
					Since all parameters are positive and the components $v_3,$ and $w_3$ are positive, it is not difficult to verify that $b>0$ with  $ m_1 = \frac{\mu \beta^*}{\sigma + \mu} $ and $ m_2 = \frac{\sigma \beta^*}{\sigma + \mu}$. Thus, it follows  from Theorem 4. 1 in \cite{Castilo_2004}  the model \eqref{eq:(2.1)} exhibits a backward bifurcation at $\mathcal{R}_0=1$ whenever $a>0$.

				\section{Proof of Theorem \ref{Theorem3.4}} \label{Appendix B}
				
				To prove the global stability of the disease-free equilibrium at $\omega=0$ and $\rho=1$, we use LaSalle Invariance Principle  \cite{$[17]$}. For this we first define a Lyapunove function  %
				\(L: \mathcal{G}\rightarrow \mathbb{R},\)
				by 
				\[L(E)=\dfrac{\mu}{3(\sigma+\mu)}A+\dfrac{\mu}{6(\sigma+\mu)}I,\]
				where 
				\[\mathcal{G} = \{ (S,V, A, I,A_1,I_1,Q, R) \in \Omega : A_1 = 0,\, I_1=0\}\subset \Omega, \quad\mbox{and}\quad E\in \mathcal{G}.\]
				Now observe that 
				\begin{align*}
					L(E_0) = 0,\qquad L(E)>0 \quad \text{for all } E\in\mathcal{G}\backslash \{E_0\}.
				\end{align*}
				Hence the function $L$ is positive definite. 
				
				We now rewrite \eqref{eq:(2.1)} in a vector form by
				\begin{align*}
					X'(t) = f(X),
				\end{align*}
				where 
				\begin{align*}
					X=&(S, V, A, I,0,0, Q, R)^T\\
					f =& (f_1 ,f_2 ,f_3 ,f_4 ,f_5 ,f_6 ,f_7 ,f_8 )^T
				\end{align*}
				with 
				\begin{align*}
					f_1 &= \Lambda- [\lambda+\sigma+\mu]S , \\
					f_2 &= \sigma S- \mu V, \\
					f_3&= \eta\lambda S  -(\theta+\gamma_1+\mu)A,\\
					f_4&= (1-\eta)\lambda S -(\epsilon+\delta+\gamma_2+\mu)I,\\
					f_5 &=  0,\\
					f_6 &=  0,\\
					f_7 &= \theta A+ \epsilon I  - (\gamma_3+\delta+\mu)Q,\\
					f_8&= \gamma_1 A + \gamma_2I+\gamma_3 Q   - \mu R.
				\end{align*}
				Let $\dot{L}$ represent the directional derivative of $L$ in the direction of  $f $. 
				Then we have
				\begin{eqnarray*}
					\dot{L}&=& \nabla L\cdot f,\\
					&=&(0,0,\dfrac{\mu}{3(\sigma+\mu)},\dfrac{\mu}{6(\sigma+\mu)},0,0,0)\cdot f ,\\
					&=&\dfrac{\mu}{3(\sigma+\mu)}\left(\eta\lambda S-k_1 A\right)+\dfrac{\mu}{6(\sigma+\mu)}\left((1-\eta)\lambda S-k_2 I\right),\\
					&\leq&\dfrac{\mu}{3(\sigma+\mu)}\left(\eta \beta(\nu A+ I) -k_1 A\right)+\dfrac{\mu}{6(\sigma+\mu)}\left((1-\eta)\beta(\nu A+I) -k_2 I\right),\\ 
					&\leq&\dfrac{\mu}{2(\sigma+\mu)}(\beta \eta \nu-k_1)A+\dfrac{\mu}{2(\sigma+\mu)}(\beta(1-\eta)-k_2)I,\\
					&\leq&\left(\dfrac{\nu B_1}{k_1}-1\right)\dfrac{\mu}{\sigma+\mu}k_1A+\left(\dfrac{B_2}{k_2}-1\right)\dfrac{\mu}{\sigma+\mu}k_2I,\\
					&\leq&\left(\dfrac{\nu B_1}{k_1}+\dfrac{B_2}{k_2}-1\right)\dfrac{\mu}{\sigma+\mu}(k_1A+k_2I),\\
					&=&( {\mathcal{R}_0}-1)\dfrac{\mu}{\sigma+\mu}\left(k_1A+k_2I\right).
				\end{eqnarray*}
				Here we used the fact that $\dfrac{S}{N}\leq1.$
				
				Thus, $\dot{L}\leq0$ on $\mathcal{G}$ whenever $ \mathcal{R}_0\leq 1$. 
				Hence, $L$ is a Lyapunov function for $E_0$ on $\mathcal{G}$. 
				Furthermore, at $E_0$ we have $\lambda = 0$. Which implies that 
				\[\dot{L}=0 \Longleftrightarrow  E=E_0.\] 
				Hence, the largest invariant set contained in $\displaystyle \mathcal{M}=\Big\{E\in \mathcal{G}: \dot{L}(E)=0\Big\}$ is $\{E_0\}$, i. e., 
				\begin{align*}
					\lim\limits_{t\to\infty}E(t) = E_0.
				\end{align*}
				Therefore, we conclude by LaSalle Invariance Principle  \cite{$[17]$} that the disease-free equilibrium  $E_0$  of  the model with $\omega=0$ and $\rho=1$ is globally asymptotically stable on $\mathcal{G}$ for $ \mathcal{R}_0 \leq1$. 
				This completes the proof of the theorem.  \hspace*{\fill} $\Box$


			\end{document}